\documentclass[]{article}
\usepackage{graphicx}
\usepackage{amsfonts}
\usepackage{amsmath}
\usepackage{amssymb}
\usepackage{fancyhdr}
\usepackage{titlesec}
\usepackage{indentfirst}
\usepackage{booktabs}
\usepackage{verbatim}
\usepackage{amsthm}
\usepackage{hyperref}
\usepackage{subcaption}
\usepackage{listings}
\usepackage[page,header]{appendix}
\usepackage{titletoc}

\numberwithin{equation}{section}
\newtheorem{theorem}{Theorem}[section]
\newtheorem{lemma}[theorem]{Lemma}

\newtheorem{proposition}[theorem]{Proposition} 

\newtheorem{remark}{Remark}[section]

\usepackage{bm}

\newcommand{\bu}{\boldsymbol{u}}
\newcommand{\bk}{\boldsymbol{k}}
\newcommand{\bv}{\boldsymbol{v}}

\newcommand{\bw}{\boldsymbol{w}}

\newcommand{\bx}{\bm{x}}

\usepackage{color}

\begin{document}

\title{On the stability of the low-rank projector-splitting integrators \\for hyperbolic and parabolic equations}

\author{Shiheng Zhang\footnote{Department of Applied Mathematics, University of Washington, Seattle, WA 98195 (shzhang3@uw.edu).} \
	\ and \ Jingwei Hu\footnote{Department of Applied Mathematics, University of Washington, Seattle, WA 98195 (hujw@uw.edu). Corresponding author.} }
	\maketitle
	
\begin{abstract}
We study the stability of low-rank projector-splitting integrators (PSI)
for linear hyperbolic and parabolic prototype equations in one spatial
and one velocity dimension. After discretization in velocity, we compare discretize-then-project (DtP) and project-then-discretize (PtD) formulations combined with finite differences in space under periodic boundary conditions. For the hyperbolic problem, while the full-tensor Lax--Friedrichs (LF) finite difference scheme with forward Euler satisfies the classical CFL condition $\nu\le 1$, we show that the corresponding PSI scheme is $L^2$-stable under a more restrictive bound $\nu\le 1/3$ in both the DtP and PtD formulations. For the parabolic problem, while the full-tensor centered finite difference scheme with backward Euler is unconditionally $L^2$-stable, we show that the corresponding PSI scheme is not. However, a hybrid PSI using backward Euler in the $K$- and $L$-steps and forward Euler in the $S$-step is unconditionally $L^2$-stable in both the DtP and PtD formulations.
\end{abstract}

{\small 

}

\section{Introduction}

In this work, we study the stability of low-rank projector-splitting integrators (PSI) for hyperbolic and parabolic equations. Our study is primarily motivated by the kinetic Vlasov-Fokker-Planck equation \cite{Villani02}, widely used to describe plasma dynamics:
\begin{equation} \label{VFP}
\partial_t f+v\cdot \nabla_x f+E(t,x)\cdot\nabla_v f=\nu_{FP} \nabla_v \cdot (T(t,x)\nabla_v f+(v-u(t,x))f), 
\end{equation}
where $f=f(t,x,v)$ is the phase-space distribution function at time $t\geq 0$, position $x\in \Omega \subset \mathbb{R}^3$, and velocity $v\in \mathbb{R}^3$. $E(t,x)$ is the electric field, which can be specified externally or determined self-consistently via Poisson's equation. The right-hand-side of (\ref{VFP}) is a Fokker-Planck operator describing charged particle collisions, where $T(t,x)$ and $u(t,x)$ denote the temperature and bulk velocity, respectively; these can also be given externally or derived from the moments of $f$. The parameter $\nu_{FP}$ is the collision frequency, which may range from $\nu_{FP} \ll 1$ (weak collisions, mildly stiff) to $\nu_{FP} \gg 1$ (strong collisions, highly stiff). 

Numerically solving equation (\ref{VFP}) is computationally expensive, as a method based on a full-tensor grid requires $O(N^6)$ complexity, where $N$ is the number of degrees of freedom in each phase-space dimension. To address this, the dynamical low-rank (DLR) method has recently emerged as a powerful on-the-fly dimension reduction technique for solving kinetic equations \cite{einkemmer2025review}. The core idea is to decompose the solution $f$ into low-rank factors that depend only on $x$ or $v$, and to evolve these factors in time. This approach reduces the computational complexity to $O(c(r)N^3)$, where $c(r)$ is a prefactor depending on the rank $r$ of the solution.

Among various DLR methods are PSI \cite{LO14}, the augmented BUG integrator \cite{CKL22}, the parallel integrator \cite{CKL24}, and the recently introduced XL integrator \cite{EHZ25}. One of the key differences between PSI and the latter three integrators is that PSI operates on a fixed-rank manifold and does not require increasing the rank during time evolution and subsequently truncating it at the end of each time step. However, a major issue with PSI is that it involves a backward-in-time substep, which can potentially lead to numerical instability and make it less suitable for stiff or dissipative problems. Moreover, there are two possible ways to implement PSI: the discretize-then-project (DtP) formulation and the project-then-discretize (PtD) formulation. In our setup, we first discretize the variable $v$, obtaining an $x$-dependent matrix system. The DtP formulation first discretizes this system in $x$ and then applies the low-rank projection, whereas the PtD formulation first applies the low-rank projection and then discretizes the resulting projected subproblems in $x$. The PtD formulation is more commonly used in the literature (e.g., \cite{EL18}). 
Nevertheless, there are situations in which DtP is more appropriate, as performing the low-rank projection may result in a non-hyperbolic system that cannot be properly handled by standard discretization methods \cite{EHW21}.

For random parabolic equations, stability of fully discrete DLR schemes was studied in \cite{kazashi2021stability}. A more closely related work is \cite{KEC23}, which studies the stability of PSI and BUG-type schemes for hyperbolic problems. Compared to this work---particularly Theorems 1 and 2 in \cite{KEC23}---our conclusions are in fact different: while the full-tensor Lax--Friedrichs (LF) finite difference scheme with forward Euler satisfies the classical CFL condition $\nu\le 1$, we show that the corresponding PSI scheme is $L^2$-stable under a more restrictive bound $\nu\le 1/3$ in both the DtP and PtD formulations. In addition, we also consider stability for parabolic equations and show that the PSI scheme, when coupled with a hybrid time discretization, can be unconditionally stable.



To simplify the discussion, we consider two prototype equations with periodic boundary conditions in the $x$-variable:
\begin{align} 
&\text{hyperbolic:} \quad \partial_t u+a(v)\partial_x u=0; \label{hyper}\\
&\text{parabolic:} \quad \partial_t u= a(v)\partial_{xx}u, \quad a(v)\geq 0. \label{parab}
\end{align}
Here $u=u(t,x,v)$ is the unknown function, and $a(v)$ denotes the wave speed in the hyperbolic equation and the diffusion coefficient in the parabolic equation (in the latter case, we also assume $a(v)\geq 0$). Equation \eqref{hyper} mimics the transport-type operator in \eqref{VFP}, while equation \eqref{parab} mimics the diffusion-type operator in \eqref{VFP}. 


We first discretize the variable $v$ in equations \eqref{hyper} and \eqref{parab}, leading to a matrix system. We then study the stability of PSI for this system when coupled with a standard finite difference discretization in $x$, carried out in either the DtP or PtD framework. 
{\bf Our main findings can be summarized as follows:}
\begin{itemize}
\item For the hyperbolic problem, the full-tensor LF scheme with forward Euler is $L^2$-stable under the classical hyperbolic CFL condition $\nu\le 1$. For the corresponding low-rank PSI schemes, single-mode calculations suggest the bound $\nu\le 1/3$, and we prove $L^2$-stability under this bound in both the DtP and PtD formulations.

\item For the parabolic problem, while the full-tensor centered finite difference scheme with backward Euler is unconditionally $L^2$-stable, we show that the corresponding PSI scheme is not. However, a hybrid PSI using backward Euler in the $K$- and $L$-steps and forward Euler in the $S$-step is unconditionally $L^2$-stable in both the DtP and PtD formulations.
\end{itemize}

The rest of this paper is organized as follows. Section~\ref{sec:prelim} introduces the $v$-discrete formulation, the structural assumptions on $A$, and the difference operators in the $x$-variable. Section~\ref{sec:hyper} studies the hyperbolic problem with the LF scheme, first in the DtP formulation and then in the PtD formulation. Section~\ref{sec:parabolic} is devoted to the parabolic problem, where we identify the hybrid DtP scheme and prove its unconditional $L^2$-stability, and then transfer the same conclusion to the PtD formulation. Finally, the paper is concluded in Section~\ref{sec:conclusion}. The appendices collect the remaining details for the hyperbolic LF analysis and the proof of the single-mode result for the DtP $\theta$-scheme in the parabolic case.

\section{Preliminaries: $v$-discrete formulation and difference matrices}
\label{sec:prelim}

We work throughout with a formulation that has already been discretized in the variable $v$. Accordingly, the unknown is an $N_v$-vector-valued function of the periodic variable $x$ and time $t$, or equivalently an $N_x\times N_v$ matrix after sampling in $x$. In this setting, the distinction between the DtP and PtD formulations concerns only the order of the low-rank projection and the discretization in the variable $x$. The purpose of this section is therefore to record the common $v$-discrete matrix formulation and to introduce the finite-difference matrices in $x$ and their basic algebraic properties, which will be used repeatedly in the following discussion.

\subsection{$v$-discrete matrix formulation}

After discretization in $v$, the hyperbolic equation \eqref{hyper} and the parabolic equation \eqref{parab} reduce to
\begin{equation}
\label{hyper-sys}
\partial_t \bu(t,x)+A\,\partial_x \bu(t,x)=0,
\qquad
\bu(t,x)\in\mathbb R^{N_v},
\end{equation}
and
\begin{equation}
\label{eq:parabolic_semi_discrete_v}
\partial_t \bu(t,x)=A\,\partial_{xx}\bu(t,x),
\end{equation}
where $A\in\mathbb R^{N_v\times N_v}$ is constant and symmetric. This setting includes, for example, nodal discretizations, for which
\[
A=\operatorname{diag}(a(v_1),\dots,a(v_{N_v})),
\]
and modal discretizations in an orthonormal basis $\{\phi_i(v)\}_{i=1}^{N_v}$, for which
\[
A_{ij}=\langle a\phi_i\phi_j\rangle_v,
\]
where $\langle\cdot\rangle_v$ denotes the integral over the domain of $v$.

For symmetric matrices $B$ and $C$, we use the Loewner partial order, writing
\[
C\succeq B
\qquad\text{or equivalently}\qquad
B\preceq C
\]
to mean that $C-B$ is positive semidefinite. In the parabolic case we additionally assume
\[
A\succeq 0,
\]
which corresponds naturally to the continuous physical condition $a(v)\ge 0$.

Accordingly, we write the spectral decomposition
\[
A=Q_A\Lambda Q_A^T,
\qquad
\Lambda=\operatorname{diag}(\lambda_1,\dots,\lambda_{N_v}),
\]
with $Q_A$ orthogonal.

\subsection{Finite-difference matrices in the periodic variable $x$}\label{sec:finite_matrices}

Although the DtP and PtD formulations discretize the periodic variable $x$ at different stages, both ultimately involve the same periodic difference operators. We therefore introduce here the corresponding difference matrices and collect the basic identities that will be used later.

Let $\{x_j\}_{j=1}^{N_x}$ be a uniform periodic grid with spacing $\Delta x$, and write
\[
\bu_j(t):=\bu(t,x_j),
\qquad j=1,\dots,N_x.
\]
To write finite-difference formulas in matrix form, we denote by $E_+$ and $E_-$ the forward and backward shift matrices, so that for any vector $\boldsymbol{w}\in\mathbb R^{N_x}$,
\[
(E_+\boldsymbol{w})_j=w_{j+1},
\qquad
(E_-\boldsymbol{w})_j=w_{j-1},
\]
with periodic indexing. We then define
\begin{equation}
\label{discretize-matrix}
M_1=E_+-E_-,
\qquad
M_2=2I-E_+-E_-.
\end{equation}
Equivalently,
\[
(M_1\boldsymbol{w})_j=w_{j+1}-w_{j-1},
\qquad
(M_2\boldsymbol{w})_j=2w_j-w_{j+1}-w_{j-1}.
\]
Hence
\[
M_1^T=-M_1,
\qquad
M_2^T=M_2,
\qquad
M_2\succeq 0,
\]
and we also have
\begin{equation}
\label{eq:M1M2}
M_1^TM_1=4M_2-M_2^2.
\end{equation}

For later reference, if one discretizes the $x$-variable in \eqref{hyper-sys} directly by the LF scheme, one obtains
\begin{equation}
\label{eq:hyper_LF_component}
\dot{\bu}_j
=
-\frac{1}{2\Delta x}A(\bu_{j+1}-\bu_{j-1})
+\frac{\lambda_{\max}}{2\Delta x}\bigl(\bu_{j-1}-2\bu_j+\bu_{j+1}\bigr),
\qquad j=1,\dots,N_x,
\end{equation}
where
\[
\lambda_{\max}:=\|A\|_2=\max_{1\le k\le N_v}|\lambda_k|,
\]
and here and below $\|\cdot\|_2$ denotes the spectral norm. Similarly, the standard second-order central difference discretization of \eqref{eq:parabolic_semi_discrete_v} leads to
\begin{equation}
\label{eq:parabolic_component}
\dot{\bu}_j
=
\frac{1}{(\Delta x)^2}A(\bu_{j-1}-2\bu_j+\bu_{j+1}),
\qquad j=1,\dots,N_x.
\end{equation}

We assemble the solution matrix
\begin{equation}
\label{UU}
U(t)=
\begin{bmatrix}
\bu_1(t)^T\\
\vdots\\
\bu_{N_x}(t)^T
\end{bmatrix}
\in\mathbb R^{N_x\times N_v}.
\end{equation}
Then \eqref{eq:hyper_LF_component} and \eqref{eq:parabolic_component} can be written as
\begin{equation}
\label{eq:hyper_LF_matrix}
\dot U
=
-\frac{1}{2\Delta x}M_1UA
-\frac{\lambda_{\max}}{2\Delta x}M_2U,
\end{equation}
and
\begin{equation}
\label{eq:parabolic_semi_discrete_x_revised3}
\dot U
=
-\frac{1}{(\Delta x)^2}M_2UA.
\end{equation}

\subsection{Matrix norms}

We equip $\mathbb R^{N_x\times N_v}$ with the Frobenius inner product
\[
\langle U,W\rangle_F:=\operatorname{tr}(U^TW),
\qquad
\|U\|_F:=\sqrt{\langle U,U\rangle_F}.
\]

The precise discrete $L^2$ norm in $(x,v)$ depends on the normalization of the $v$-discretization. In the orthonormal modal setting,
\[
\|U\|_{L_h^2(x,v)}^2=\Delta x\,\|U\|_F^2.
\]
In a nodal discretization with quadrature weights $(\omega_\ell)_{\ell=1}^{N_v}$,
\[
\|U\|_{L_h^2(x,v)}^2
=
\Delta x\sum_{j=1}^{N_x}\sum_{\ell=1}^{N_v}\omega_\ell |U_{j\ell}|^2.
\]
Thus, after fixing the discretization in $v$, the discrete $L^2$ norm differs from $\|\cdot\|_F$ only by a fixed normalization (or, more generally, fixed positive weights). In particular, all stability statements below can be equivalently formulated in terms of $\|\cdot\|_F$.

We also use the discrete $H_x^1$-seminorm
\[
|U|_{1,x}^2:=\langle U,M_2U\rangle_F.
\]
By periodicity,
\[
|U|_{1,x}^2
=
\sum_{j=1}^{N_x} |\bu_{j+1}-\bu_j|^2.
\]
This differs from the standard grid-scaled discrete $H_x^1$ seminorm only by a fixed normalization factor, which is immaterial for the estimates below.

\section{Hyperbolic problem}
\label{sec:hyper}

We now turn to the hyperbolic semi-discretization \eqref{eq:hyper_LF_matrix}. For the full-tensor LF scheme coupled with forward Euler, the classical CFL condition is $\nu\le 1$. Our goal in this section is to determine how this stability picture changes for the corresponding low-rank PSI schemes. We consider the DtP and PtD formulations separately. In each case, we first write the corresponding low-rank scheme, then examine a single Fourier mode to identify the bound suggested by the mode amplification, and finally prove an $L^2$-stability result for general low-rank states.

\subsection{Full-tensor scheme}
\label{subsec:hyper_full_lf}

Applying forward Euler to \eqref{eq:hyper_LF_matrix} gives
\begin{equation}
\label{eq:hyper_full_fe}
U^{n+1}
=
U^n-\frac{\Delta t}{2\Delta x}M_1U^nA-\frac{\lambda_{\max}\Delta t}{2\Delta x}M_2U^n.
\end{equation}
We introduce the hyperbolic CFL number
\begin{equation}
\label{eq:hyper_cfl}
\nu:=\lambda_{\max}\frac{\Delta t}{\Delta x}>0.
\end{equation}
Throughout this section the periodic indexing convention introduced in
Section~\ref{sec:finite_matrices} is in effect; in particular, the shift matrices $E_\pm$ and
the difference matrices $M_1$, $M_2$ are all defined with respect to the
periodic grid $\{x_j\}_{j=1}^{N_x}$.
The following standard estimate provides the baseline stability condition
for the low-rank analysis below.

\begin{theorem}
\label{thm:hyper_full_lf}
If $\nu\le 1$, then the full-tensor scheme \eqref{eq:hyper_full_fe} satisfies
\[
\|U^{n+1}\|_F\le \|U^n\|_F,
\qquad \text{for all } n\ge 0.
\]
\end{theorem}

\begin{proof}
Let $A=Q_A\Lambda Q_A^T$ be the spectral decomposition from Section~\ref{sec:prelim}, with
\[
\Lambda=\operatorname{diag}(\lambda_1,\dots,\lambda_{N_v}),
\]
and define
\[
W^n:=U^nQ_A.
\]
Since $Q_A$ is orthogonal, $\|W^n\|_F=\|U^n\|_F$. Multiplying \eqref{eq:hyper_full_fe} on the right by $Q_A$, we obtain
\[
W^{n+1}
=
W^n-\frac{\Delta t}{2\Delta x}M_1W^n\Lambda-\frac{\nu}{2}M_2W^n.
\]
Thus each column $\bw_j^n\in\mathbb R^{N_x}$ of $W^n$ satisfies
\[
\bw_j^{n+1}
=
\bw_j^n-\frac{\nu_j}{2}M_1\bw_j^n-\frac{\nu}{2}M_2\bw_j^n,
\qquad
\nu_j:=\lambda_j\frac{\Delta t}{\Delta x}.
\]
Using $M_1=E_+-E_-$ and $M_2=2I-E_+-E_-$, this becomes
\[
\bw_j^{n+1}
=
\frac{\nu-\nu_j}{2}E_+\bw_j^n
+
(1-\nu)\bw_j^n
+
\frac{\nu+\nu_j}{2}E_-\bw_j^n.
\]
If $\nu\le 1$, then $|\nu_j|\le \nu$, so the three coefficients are nonnegative and sum to one.
Since $E_+$ and $E_-$ are orthogonal matrices,
\[
\|\bw_j^{n+1}\|_2
\le
\left(
\frac{\nu-\nu_j}{2}+(1-\nu)+\frac{\nu+\nu_j}{2}
\right)\|\bw_j^n\|_2
=
\|\bw_j^n\|_2.
\]
Summing over $j$ yields
\[
\|U^{n+1}\|_F=\|W^{n+1}\|_F\le \|W^n\|_F=\|U^n\|_F.
\]
\end{proof}


\subsection{DtP formulation}
\label{subsec:hyper_dtp}

We begin with the DtP formulation. Starting from the matrix equation \eqref{eq:hyper_LF_matrix}, we apply the PSI scheme using the Lie-Trotter splitting and forward Euler time stepping in each substep.

\subsubsection{The DtP--PSI scheme}
\label{subsec:dtp_lf_scheme}

We consider the low-rank approximation
\begin{equation}
\label{eq:hyper_low_rank}
U(t)=X(t)S(t)V(t)^T,
\qquad
X\in\mathbb R^{N_x\times r},\quad
S\in\mathbb R^{r\times r},\quad
V\in\mathbb R^{N_v\times r},
\end{equation}
with $X^TX=I_r$ and $V^TV=I_r$. For the matrix equation \eqref{eq:hyper_LF_matrix}, define
\[
F_{\rm LF}(U):=-\frac{1}{2\Delta x}M_1UA-\frac{\lambda_{\max}}{2\Delta x}M_2U.
\]
The DtP--PSI scheme from $t^n$ to $t^{n+1}$, with forward Euler in each substep, is
\begin{itemize}
\item K-step: let $K^n=X^nS^n$ and compute
\begin{equation}
\label{eq:dtp_lf_K}
K^{n+1}=K^n+\Delta t F_{\rm LF}(K^n(V^n)^T)V^n.
\end{equation}
Then take a QR factorization $K^{n+1}=X^{n+1}S^{(1)}$.

\item S-step:
\begin{equation}
\label{eq:dtp_lf_S}
S^{(2)}=S^{(1)}-\Delta t (X^{n+1})^TF_{\rm LF}(X^{n+1}S^{(1)}(V^n)^T)V^n.
\end{equation}

\item L-step: let $L^n=S^{(2)}(V^n)^T$ and compute
\begin{equation}
\label{eq:dtp_lf_L}
L^{n+1}=L^n+\Delta t (X^{n+1})^TF_{\rm LF}(X^{n+1}L^n).
\end{equation}
Then take a QR factorization $(L^{n+1})^T=V^{n+1}(S^{n+1})^T$.
\end{itemize}

It is convenient to introduce the intermediate solutions
\[
U^{(1)}:=X^{n+1}S^{(1)}(V^n)^T,
\qquad
U^{(2)}:=X^{n+1}S^{(2)}(V^n)^T,
\]
and the orthogonal projectors
\[
P_{V^n}:=V^n(V^n)^T,
\qquad
P_{X^{n+1}}:=X^{n+1}(X^{n+1})^T.
\]
Then \eqref{eq:dtp_lf_K}--\eqref{eq:dtp_lf_L} are equivalent to
\begin{align}
U^{(1)}
&=
U^n-\frac{\Delta t}{2\Delta x}\bigl(M_1U^nA+\lambda_{\max}M_2U^n\bigr)P_{V^n},
\label{eq:dtp_lf_full_1}
\\
U^{(2)}
&=
U^{(1)}+\frac{\Delta t}{2\Delta x}P_{X^{n+1}}\bigl(M_1U^{(1)}A+\lambda_{\max}M_2U^{(1)}\bigr)P_{V^n},
\label{eq:dtp_lf_full_2}
\\
U^{n+1}
&=
U^{(2)}-\frac{\Delta t}{2\Delta x}P_{X^{n+1}}\bigl(M_1U^{(2)}A+\lambda_{\max}M_2U^{(2)}\bigr).
\label{eq:dtp_lf_full_3}
\end{align}

We first examine the action of this scheme on a single Fourier mode and derive the corresponding CFL condition.

\subsubsection{Single-mode analysis}
\label{subsec:dtp_single_mode}

For $m=0,\dots,N_x-1$, define the complex Fourier mode
\[
(\bx_m)_j:=\frac{1}{\sqrt{N_x}}e^{ij\alpha_m},
\qquad
\alpha_m:=\frac{2\pi m}{N_x},
\qquad j=1,\dots,N_x.
\]
Then
\[
E_\pm\bx_m=e^{\pm i\alpha_m}\bx_m,
\qquad
M_1\bx_m=2iz_m\bx_m,
\qquad
M_2\bx_m=2y_m\bx_m,
\]
where
\[
y_m:=1-\cos\alpha_m,
\qquad
z_m:=\sin\alpha_m,
\qquad
z_m^2=y_m(2-y_m).
\]

\begin{proposition}
\label{prop:dtp_single_mode}
Assume that $U^n=\bx_m\bv_k^T$, where $A\bv_k=\lambda_k\bv_k$ and $\|\bv_k\|_2=1$. Then one DtP--PSI step gives
\[
U^{n+1}=g^{\rm DtP}_{m,k}U^n,
\qquad
g^{\rm DtP}_{m,k}=p_{m,k}^2(2-p_{m,k}),
\]
with
\[
p_{m,k}:=1-\nu y_m-i\nu_k z_m,
\qquad
\nu_k:=\lambda_k\frac{\Delta t}{\Delta x}.
\]
Consequently,
\[
|g^{\rm DtP}_{m,k}|^2
=
\bigl((1-\nu y_m)^2+\nu_k^2 z_m^2\bigr)^2
\bigl((1+\nu y_m)^2+\nu_k^2 z_m^2\bigr).
\]
If, in addition, $|\lambda_k|=\lambda_{\max}$, then $|\nu_k|=\nu$ and
\[
|g^{\rm DtP}_{m,k}|^2
=
h_{\rm DtP}(y_m,\nu)
:=
\bigl[1+2y_m\nu(\nu-1)\bigr]^2
\bigl[1+2y_m\nu(\nu+1)\bigr].
\]
Moreover, $h_{\rm DtP}(y,\nu)\le 1$ for all $y\in[0,2]$ when $0\le \nu\le \frac13$, whereas $h_{\rm DtP}(y,\nu)>1$ for all sufficiently small $y>0$ when $\nu>\frac13$.
\end{proposition}

The proof is given in Appendix~\ref{app:dtp_single_mode_calc}. 

Proposition~\ref{prop:dtp_single_mode} shows that the single-mode analysis leads to the bound $\nu\le \frac13$. Since rank-one single-mode states form a special class of general low-rank states, this bound is already necessary for any general $L^2$-stability result. However, unlike the classical full-tensor scheme, this single-mode calculation does not by itself imply $L^2$-stability for arbitrary low-rank states. For the full-tensor scheme, the one-step map is a fixed linear operator, so the Fourier modes in $x$ evolve independently and the result extends by superposition. This is no longer true for the PSI scheme: the projectors $P_{X^{n+1}}$ and $P_{V^n}$, as well as the QR factorizations, depend on the current iterate. As a result, the one-step map is not a fixed linear amplification operator, and the projection and reorthogonalization steps couple the modes. In the following, we show that the same CFL condition derived from the single-mode analysis also carries over to the general low-rank states.

\subsubsection{$L^2$-stability for general low-rank states}
\label{subsec:dtp_general}

\begin{theorem}
\label{thm:dtp_general}
Consider the DtP--PSI scheme \eqref{eq:dtp_lf_K}--\eqref{eq:dtp_lf_L}. If $\nu\le \frac13$, then
\[
\|U^{n+1}\|_F\le \|U^n\|_F,
\qquad \text{for all } n\ge 0.
\]
\end{theorem}

The proof starts from a mixed $L^2$--$H^1_x$ estimate produced by the K-step and then controls the S-step and L-step in reduced form. We first state a lemma that will be used in the proof of the theorem.

\begin{lemma}
\label{lem:dtp_metric}
Let $\mathcal E:\mathbb R^{r\times r}\to\mathbb R^{r\times r}$ be linear and satisfy
\[
\|\mathcal E\|_{F\to F}\le 1,
\]
where $\|\cdot\|_{F\to F}$ is the operator norm induced by the Frobenius norm. Define
\[
\mathcal H:=2\mathcal I-(\mathcal E+\mathcal E^*),
\qquad
\mathcal W:=\mathcal I-\nu(1-\nu)\mathcal H,
\qquad
\mathcal Q:=\mathcal I+\nu(\mathcal I-\mathcal E),
\]
where $\mathcal E^*$ is the adjoint with respect to the Frobenius inner product and $\mathcal I$ is the identity operator on $\mathbb R^{r\times r}$. If $\nu\le \frac13$, then $\mathcal W\succ 0$ and
\[
\mathcal Q^*\mathcal W\mathcal Q
\preceq
\mathcal I+\nu(1-\nu)\mathcal H.
\]
\end{lemma}

The proof is given in Appendix~\ref{app:dtp_metric_proof}.

\begin{proof}
Write the discrete $H^1_x$-seminorm as
\[
|U|_{1,x}^2:=\langle U,M_2U\rangle_F.
\]

\paragraph{Step 1.}
Define
\[
G_0:=I-\frac{\nu}{2}M_2,
\qquad
G_1:=\frac{\Delta t}{2\Delta x}M_1,
\qquad
P_{V^n}:=V^n(V^n)^T.
\]
Then
\[
G_0^T=G_0,
\qquad
G_1^T=-G_1,
\qquad
G_0G_1=G_1G_0.
\]
Using \eqref{eq:dtp_lf_full_1} and $U^nP_{V^n}=U^n$, the K-step can be written as
\begin{equation}
\label{eq:dtp_k_rewrite}
U^{(1)}=(G_0U^n-G_1U^nA)P_{V^n}.
\end{equation}
Since $P_{V^n}$ is an orthogonal projector,
\[
\|U^{(1)}\|_F^2
\le
\|G_0U^n-G_1U^nA\|_F^2.
\]
Expanding the Frobenius norm gives
\[
\|G_0U^n-G_1U^nA\|_F^2
=
\|G_0U^n\|_F^2+\|G_1U^nA\|_F^2
-2\langle G_0U^n,G_1U^nA\rangle_F.
\]
Moreover,
\[
\langle G_0U^n,G_1U^nA\rangle_F
=
\operatorname{tr}\!\bigl((U^n)^TG_0^TG_1U^nA\bigr)=0,
\]
because $G_0^TG_1$ is skew-symmetric, and hence
\(
(U^n)^TG_0^TG_1U^n
\)
is also skew-symmetric, while $A$ is symmetric.
Therefore
\[
\|G_0U^n-G_1U^nA\|_F^2
=
\|G_0U^n\|_F^2+\|G_1U^nA\|_F^2.
\]
Using $A^2\preceq \lambda_{\max}^2I$ and
\(
(U^n)^TG_1^TG_1U^n\succeq 0,
\)
we obtain
\[
\|G_1U^nA\|_F^2
=
\operatorname{tr}\!\bigl((U^n)^TG_1^TG_1U^nA^2\bigr)
\le
\lambda_{\max}^2\|G_1U^n\|_F^2.
\]
Hence
\[
\|U^{(1)}\|_F^2
\le
\bigl\langle U^n,\bigl(G_0^TG_0+\lambda_{\max}^2G_1^TG_1\bigr)U^n\bigr\rangle_F.
\]
Since
\[
G_0^TG_0
=
\Bigl(I-\frac{\nu}{2}M_2\Bigr)^2
=
I-\nu M_2+\frac{\nu^2}{4}M_2^2
\]
and
\[
\lambda_{\max}^2G_1^TG_1
=
\lambda_{\max}^2\Bigl(\frac{\Delta t}{2\Delta x}\Bigr)^2M_1^TM_1
=
\frac{\nu^2}{4}(4M_2-M_2^2)
=
\nu^2M_2-\frac{\nu^2}{4}M_2^2,
\]
identity \eqref{eq:M1M2} yields
\[
G_0^TG_0+\lambda_{\max}^2G_1^TG_1
=
I-\nu(1-\nu)M_2.
\]
Therefore
\begin{equation}
\label{eq:dtp_after_k}
\|U^{(1)}\|_F^2
\le
\|U^n\|_F^2-\nu(1-\nu)|U^n|_{1,x}^2.
\end{equation}

Next we consider $M_2^{1/2}U^n$. Since $M_1$ and $M_2$ are circulant, they are simultaneously diagonalized by the discrete Fourier matrix, so $M_2^{1/2}$ commutes with $M_1$, $G_0$, and $G_1$. Multiplying \eqref{eq:dtp_k_rewrite} by $M_2^{1/2}$ from the left gives
\[
M_2^{1/2}U^{(1)}=(G_0M_2^{1/2}U^n-G_1M_2^{1/2}U^nA)P_{V^n}.
\]
Repeating the previous estimate yields
\[
|U^{(1)}|_{1,x}^2
=
\|M_2^{1/2}U^{(1)}\|_F^2
\le
\|M_2^{1/2}U^n\|_F^2-\nu(1-\nu)|M_2^{1/2}U^n|_{1,x}^2
\le
|U^n|_{1,x}^2,
\]
because $M_2\succeq 0$. Hence
\begin{equation}
\label{eq:dtp_mixed_energy}
\|U^{(1)}\|_F^2+\nu(1-\nu)|U^{(1)}|_{1,x}^2
\le
\|U^n\|_F^2.
\end{equation}

\paragraph{Step 2.}
After the K-step, write
\[
U^{(1)}=X^{n+1}S^{(1)}(V^n)^T,
\qquad
\widetilde A:=(V^n)^TAV^n,
\qquad
\widehat A:=\frac{1}{\lambda_{\max}}\widetilde A.
\]
Since $A$ is symmetric and $\lambda_{\max}=\|A\|_2$, the reduced matrix $\widetilde A$ is symmetric and satisfies
\[
\widetilde A^2\preceq \lambda_{\max}^2I,
\qquad
\widehat A^T=\widehat A,
\qquad
\widehat A^2\preceq I.
\]
Define the reduced matrices
\[
M_{1,X}:=(X^{n+1})^TM_1X^{n+1},
\qquad
M_{2,X}:=(X^{n+1})^TM_2X^{n+1}.
\]
Then
\[
M_{1,X}^T=-M_{1,X},
\qquad
0\preceq M_{2,X}\preceq 4I.
\]
Moreover,
\begin{equation}
\label{eq:dtp_reduced_h1}
|U^{(1)}|_{1,x}^2
=
\operatorname{tr}\!\bigl((S^{(1)})^TM_{2,X}S^{(1)}\bigr).
\end{equation}
Indeed,
\[
|U^{(1)}|_{1,x}^2
=
\bigl\langle X^{n+1}S^{(1)}(V^n)^T,\,
M_2X^{n+1}S^{(1)}(V^n)^T\bigr\rangle_F
=
\operatorname{tr}\!\bigl((S^{(1)})^TM_{2,X}S^{(1)}\bigr).
\]

Let $\mathcal I$ denote the identity operator on any space under consideration, and take all adjoints with respect to the Frobenius inner product. Define the full LF transport operator $\mathcal T:\mathbb R^{N_x\times N_v}\to\mathbb R^{N_x\times N_v}$ by
\[
\mathcal T(Y)
:=
\Bigl(I-\frac12M_2\Bigr)Y-\frac{1}{2\lambda_{\max}}M_1YA
=
\frac12E_+Y\Bigl(I-\frac{A}{\lambda_{\max}}\Bigr)
+
\frac12E_-Y\Bigl(I+\frac{A}{\lambda_{\max}}\Bigr).
\]
Using $E_+^T=E_-$, $E_-^T=E_+$, $E_+E_-=E_-E_+=I$, and $A=A^T$, we have
\[
\mathcal T^*(Y)
=
\frac12E_-Y\Bigl(I-\frac{A}{\lambda_{\max}}\Bigr)
+
\frac12E_+Y\Bigl(I+\frac{A}{\lambda_{\max}}\Bigr).
\]
A direct expansion gives
\[
\mathcal T^*\mathcal T(Y)
=
\frac12Y\Bigl(I+\frac{A^2}{\lambda_{\max}^2}\Bigr)
+
\frac14(E_+^2+E_-^2)Y\Bigl(I-\frac{A^2}{\lambda_{\max}^2}\Bigr).
\]
Hence, using $2I-E_+^2-E_-^2=M_1^TM_1$, we obtain
\begin{equation}
\label{eq:dtp_transport_contraction}
\|Y\|_F^2-\|\mathcal T(Y)\|_F^2
=
\frac14\operatorname{tr}\!\Bigl(Y^TM_1^TM_1Y\Bigl(I-\frac{A^2}{\lambda_{\max}^2}\Bigr)\Bigr)
\ge 0
\end{equation}
for all $Y\in\mathbb R^{N_x\times N_v}$. Therefore $\|\mathcal T\|_{F\to F}\le 1$.

Compressing on the left by $X^{n+1}$, define $\mathcal T_X:\mathbb R^{r\times N_v}\to\mathbb R^{r\times N_v}$ by
\[
\mathcal T_X(L):=(X^{n+1})^T\mathcal T(X^{n+1}L)
=
\Bigl(I-\frac12M_{2,X}\Bigr)L-\frac{1}{2\lambda_{\max}}M_{1,X}LA.
\]
Since $X^{n+1}$ has orthonormal columns,
\[
\|(X^{n+1})^TY\|_F\le \|Y\|_F,
\qquad
\|X^{n+1}L\|_F=\|L\|_F.
\]
Therefore,
\[
\|\mathcal T_X(L)\|_F
\le
\|\mathcal T(X^{n+1}L)\|_F
\le
\|X^{n+1}L\|_F
=
\|L\|_F,
\]
so
\begin{equation}
\label{eq:dtp_reduced_transport_contraction}
\|\mathcal T_X\|_{F\to F}\le 1.
\end{equation}

Using \eqref{eq:dtp_lf_L} and $L^n=S^{(2)}(V^n)^T$, the L-step becomes
\[
L^{n+1}=L^n-\frac{\nu}{2\lambda_{\max}}M_{1,X}L^nA-\frac{\nu}{2}M_{2,X}L^n=:\mathcal B(L^n),
\qquad
\mathcal B:=(1-\nu)\mathcal I+\nu\mathcal T_X.
\]
Set
\begin{equation}
\label{eq:dtp_WX}
W_X:=I-\nu(1-\nu)M_{2,X}.
\end{equation}
Because $A^T=A$ and $M_{1,X}^T=-M_{1,X}$, for every $L\in\mathbb R^{r\times N_v}$ we have
\[
\langle L,M_{1,X}LA\rangle_F=\operatorname{tr}(L^TM_{1,X}LA)=0,
\]
and therefore
\begin{equation}
\label{eq:dtp_lstep_metric}
\|\mathcal B(L)\|_F^2+\nu^2\bigl(\|L\|_F^2-\|\mathcal T_X(L)\|_F^2\bigr)
=
\langle L,W_XL\rangle_F.
\end{equation}
Using \eqref{eq:dtp_reduced_transport_contraction}, we obtain from \eqref{eq:dtp_lstep_metric} that
\begin{equation}
\label{eq:dtp_lstep_bound}
\|U^{n+1}\|_F^2
=
\|L^{n+1}\|_F^2
\le
\operatorname{tr}\!\bigl((S^{(2)})^TW_XS^{(2)}\bigr).
\end{equation}

To control the S-step, define $\mathcal E:\mathbb R^{r\times r}\to\mathbb R^{r\times r}$ by
\[
\mathcal E(Z)
:=
\Bigl(I-\frac12M_{2,X}\Bigr)Z-\frac12M_{1,X}Z\widehat A.
\]
Then
\[
\mathcal E(Z)=\mathcal T_X\bigl(Z(V^n)^T\bigr)V^n.
\]
Since $V^n$ has orthonormal columns,
\[
\|WV^n\|_F\le \|W\|_F
\qquad\text{for all }W\in\mathbb R^{r\times N_v}.
\]
Therefore,
\[
\|\mathcal E(Z)\|_F
=
\|\mathcal T_X\bigl(Z(V^n)^T\bigr)V^n\|_F
\le
\|\mathcal T_X\bigl(Z(V^n)^T\bigr)\|_F
\le
\|Z(V^n)^T\|_F
=
\|Z\|_F,
\]
so
\begin{equation}
\label{eq:dtp_S_operator_contraction}
\|\mathcal E\|_{F\to F}\le 1.
\end{equation}
Using \eqref{eq:dtp_lf_S}, the S-step becomes
\begin{equation}
\label{eq:dtp_S_as_Q}
S^{(2)}=S^{(1)}+\frac{\nu}{2}M_{1,X}S^{(1)}\widehat A+\frac{\nu}{2}M_{2,X}S^{(1)}=\mathcal Q(S^{(1)}),
\qquad
\mathcal Q:=\mathcal I+\nu(\mathcal I-\mathcal E).
\end{equation}

\paragraph{Step 3.}
Define
\[
\mathcal H:=2\mathcal I-(\mathcal E+\mathcal E^*),
\qquad
\mathcal W:=\mathcal I-\nu(1-\nu)\mathcal H.
\]
Since
\[
\mathcal E^*(Z)=\Bigl(I-\frac12M_{2,X}\Bigr)Z+\frac12M_{1,X}Z\widehat A,
\]
we have
\[
\mathcal H(Z)=M_{2,X}Z,
\qquad
\mathcal W(Z)=W_XZ.
\]
Combining \eqref{eq:dtp_lstep_bound} and \eqref{eq:dtp_S_as_Q}, we obtain
\[
\|U^{n+1}\|_F^2
\le
\langle \mathcal Q(S^{(1)}),\mathcal W\mathcal Q(S^{(1)})\rangle_F.
\]
Apply Lemma~\ref{lem:dtp_metric} with the operator $\mathcal E$. By \eqref{eq:dtp_S_operator_contraction},
\[
\mathcal Q^*\mathcal W\mathcal Q
\preceq
\mathcal I+\nu(1-\nu)\mathcal H.
\]
Therefore, using \eqref{eq:dtp_reduced_h1},
\[
\|U^{n+1}\|_F^2
\le
\|S^{(1)}\|_F^2+\nu(1-\nu)\operatorname{tr}\!\bigl((S^{(1)})^TM_{2,X}S^{(1)}\bigr)
=
\|U^{(1)}\|_F^2+\nu(1-\nu)|U^{(1)}|_{1,x}^2.
\]
Finally, \eqref{eq:dtp_mixed_energy} yields
\[
\|U^{n+1}\|_F^2\le \|U^n\|_F^2.
\]
This proves the $L^2$-stability of the DtP--PSI scheme under $\nu\le \frac13$.
\end{proof}

\begin{remark}[A second-order DtP variant: single-mode calculation only]
Consider the PSI scheme using the Strang splitting
\[
\text{half K-step} \;\to\; \text{half S-step} \;\to\; \text{full L-step} \;\to\; \text{half S-step} \;\to\; \text{half K-step},
\]
and apply the second-order strong-stability-preserving Runge--Kutta method (SSP-RK2) \cite{GKS11} in each substep. For a scalar one-step factor $q$, the associated SSP-RK2 factor is
\[
R_{\mathrm{SSP2}}(q):=\frac12(1+q^2).
\]

Assume that $U^n=\bx_m\bv_k^T$ as in Proposition~\ref{prop:dtp_single_mode}, and keep the notation there. For the half-step quantities, let
\[
p_{m,k}^{(\frac12)}:=1-\frac{\nu}{2}y_m-i\frac{\nu_k}{2}z_m.
\]
Then one such step gives
\[
U^{n+1}=g^{\rm DtP,Strang}_{m,k}U^n,
\]
with amplification factor
\[
g^{\rm DtP,Strang}_{m,k}
=
R_{\mathrm{SSP2}}\!\bigl(p_{m,k}^{(\frac12)}\bigr)\,
R_{\mathrm{SSP2}}\!\bigl(2-p_{m,k}^{(\frac12)}\bigr)\,
R_{\mathrm{SSP2}}(p_{m,k})\,
R_{\mathrm{SSP2}}\!\bigl(2-p_{m,k}^{(\frac12)}\bigr)\,
R_{\mathrm{SSP2}}\!\bigl(p_{m,k}^{(\frac12)}\bigr).
\]
Equivalently,
\[
|g^{\rm DtP,Strang}_{m,k}|^2
=
\bigl|R_{\mathrm{SSP2}}(p_{m,k}^{(\frac12)})\bigr|^4
\bigl|R_{\mathrm{SSP2}}(2-p_{m,k}^{(\frac12)})\bigr|^4
\bigl|R_{\mathrm{SSP2}}(p_{m,k})\bigr|^2.
\]

If, in addition, $|\lambda_k|=\lambda_{\max}$, then $|\nu_k|=\nu$, and the above expression becomes a scalar function of $y_m\in[0,2]$ and $\nu$. A direct numerical scan of this function indicates that
\[
|g^{\rm DtP,Strang}_{m,k}|\le 1
\qquad\text{for all } y_m\in[0,2]
\]
whenever
\[
\nu \lesssim 0.866.
\]
Although the above analysis only works for the single mode, it still provides a useful necessary CFL condition. We do not claim a general $L^2$-stability result for arbitrary low-rank states.
\end{remark}

\subsection{PtD formulation}
\label{subsec:hyper_ptd}

We next consider the PtD formulation. We first apply the PSI to the $x$-continuous system \eqref{hyper-sys} and then discretize the resulting subproblems in $x$ by the LF scheme. 

\subsubsection{The PtD--PSI scheme}
\label{subsec:ptd_lf_scheme}

Since $A$ is symmetric, it is convenient to work with the transposed form of \eqref{hyper-sys},
\[
\partial_t\bu^T(t,x)=-\partial_x\bu^T(t,x)A.
\]
We consider the rank-$r$ approximation
\[
\bu^T(t,x)=\bx^T(t,x)S(t)V(t)^T,
\]
where $\bx(t,x)\in\mathbb R^r$, $S(t)\in\mathbb R^{r\times r}$, and $V(t)\in\mathbb R^{N_v\times r}$ satisfy
\[
\langle \bx\bx^T\rangle_x=I_r,
\qquad
V^TV=I_r.
\]
Let
\[
\bk^T(t,x):=\bx^T(t,x)S(t),
\qquad
L(t):=S(t)V(t)^T.
\]
The corresponding projector-splitting subproblems are
\begin{align}
\partial_t\bk^T &= -\partial_x\bk^T\widetilde A,
\qquad
\widetilde A:=V^TAV,
\label{eq:ptd_cont_K}
\\
\dot S &= \langle \bx\,\partial_x\bx^T\rangle_x S\widetilde A,
\label{eq:ptd_cont_S}
\\
\dot L &= -\langle \bx\,\partial_x\bx^T\rangle_x LA.
\label{eq:ptd_cont_L}
\end{align}

Among the above subproblems, only the $K$-step is a hyperbolic equation in $x$. In the $S$-step and $L$-step, the dependence on $x$ enters only through the reduced matrix $\langle \bx\,\partial_x\bx^T\rangle_x$, rather than through a full transport equation. Under periodic boundary conditions, integration by parts shows that this reduced coefficient is skew-symmetric. Its centered-difference approximation preserves the same skew-symmetry at the discrete level. We therefore discretize the $K$-step by the LF scheme and approximate $\langle \bx\,\partial_x\bx^T\rangle_x$ by the centered finite-difference operator on the periodic grid. 

Let $\{x_j\}_{j=1}^{N_x}$ be the periodic grid from Section~\ref{sec:prelim}, and write
\[
\bx_j(t):=\bx(t,x_j).
\]
Define
\[
X(t):=\sqrt{\Delta x}
\begin{bmatrix}
\bx_1(t)^T\\
\vdots\\
\bx_{N_x}(t)^T
\end{bmatrix}
\in\mathbb R^{N_x\times r},
\qquad
K(t):=X(t)S(t).
\]
Then $X(t)^TX(t)=I_r$. Applying the LF discretization to \eqref{eq:ptd_cont_K}, replacing $\langle \bx\,\partial_x\bx^T\rangle_x$ by $\frac{1}{2\Delta x}X^TM_1X$, and using forward Euler in time, we obtain 
\begin{equation}
\label{eq:ptd_lf_K}
K^{n+1}
=
K^n-\frac{\Delta t}{2\Delta x}M_1K^n\widetilde A-\frac{\nu}{2}M_2K^n,
\qquad
\widetilde A:=(V^n)^TAV^n,
\end{equation}
followed by the QR factorization
\[
K^{n+1}=X^{n+1}S^{(1)}.
\]
The S-step becomes
\begin{equation}
\label{eq:ptd_lf_S}
S^{(2)}
=
S^{(1)}+\frac{\Delta t}{2\Delta x}(X^{n+1})^TM_1X^{n+1}S^{(1)}\widetilde A,
\end{equation}
and, with
\[
L^n:=S^{(2)}(V^n)^T,
\]
the L-step is
\begin{equation}
\label{eq:ptd_lf_L}
L^{n+1}
=
L^n-\frac{\Delta t}{2\Delta x}(X^{n+1})^TM_1X^{n+1}L^nA,
\end{equation}
followed by the QR factorization
\[
(L^{n+1})^T=V^{n+1}(S^{n+1})^T.
\]

Using again the intermediate solutions
\[
U^{(1)}:=X^{n+1}S^{(1)}(V^n)^T,
\qquad
U^{(2)}:=X^{n+1}S^{(2)}(V^n)^T,
\]
we can rewrite \eqref{eq:ptd_lf_K}--\eqref{eq:ptd_lf_L} as
\begin{align}
U^{(1)}
&=
U^n-\frac{\Delta t}{2\Delta x}
\Bigl(M_1(U^nV^n)\widetilde A+\lambda_{\max}M_2(U^nV^n)\Bigr)(V^n)^T,
\label{eq:ptd_lf_full_1}
\\
U^{(2)}
&=
U^{(1)}+\frac{\Delta t}{2\Delta x}
X^{n+1}(X^{n+1})^TM_1(U^{(1)}V^n)\widetilde A(V^n)^T,
\label{eq:ptd_lf_full_2}
\\
U^{n+1}
&=
U^{(2)}-\frac{\Delta t}{2\Delta x}
X^{n+1}(X^{n+1})^T(M_1U^{(2)}A).
\label{eq:ptd_lf_full_3}
\end{align}

We next examine the action of this scheme on a single Fourier mode.

\subsubsection{Single-mode analysis}
\label{subsec:ptd_single_mode}

We use the Fourier modes $\bx_m$, together with the scalars $y_m$ and $z_m$, introduced in Section~\ref{subsec:dtp_single_mode}.

\begin{proposition}
\label{prop:ptd_single_mode}
Assume that $U^n=\bx_m\bv_k^T$, where $A\bv_k=\lambda_k\bv_k$ and $\|\bv_k\|_2=1$. Then one PtD--PSI step gives
\[
U^{n+1}=g^{\rm PtD}_{m,k}U^n,
\qquad
g^{\rm PtD}_{m,k}=p_{m,k}(1+i\nu_k z_m)(1-i\nu_k z_m),
\]
with
\[
p_{m,k}:=1-\nu y_m-i\nu_k z_m,
\qquad
\nu_k:=\lambda_k\frac{\Delta t}{\Delta x}.
\]
Consequently,
\[
|g^{\rm PtD}_{m,k}|^2
=
\bigl((1-\nu y_m)^2+\nu_k^2z_m^2\bigr)
\bigl(1+\nu_k^2 z_m^2\bigr)^2.
\]
If, in addition, $|\lambda_k|=\lambda_{\max}$, then $|\nu_k|=\nu$ and
\[
|g^{\rm PtD}_{m,k}|^2
=
h_{\rm PtD}(y_m,\nu)
:=
\bigl[1+2y_m\nu(\nu-1)\bigr]
\bigl[1+\nu^2 y_m(2-y_m)\bigr]^2.
\]
Moreover, $h_{\rm PtD}(y,\nu)\le 1$ for all $y\in[0,2]$ when $0\le \nu\le \frac13$, whereas $h_{\rm PtD}(y,\nu)>1$ for all sufficiently small $y>0$ when $\nu>\frac13$.
\end{proposition}

The proof is given in Appendix~\ref{app:ptd_single_mode_calc}. 

Proposition~\ref{prop:ptd_single_mode} shows that the single-mode analysis again leads to the bound $\nu\le \frac13$. This bound is necessary for any general $L^2$-stability result, since rank-one single-mode states form a special class of general low-rank states. However, the same caveat applies here: the single-mode calculation does not by itself imply $L^2$-stability for arbitrary low-rank states. In the following, we show that the same CFL condition derived from the single-mode analysis also carries over to the general low-rank states.

\subsubsection{$L^2$-stability for general low-rank states}
\label{subsec:ptd_general}

\begin{theorem}
\label{thm:ptd_general}
Consider the PtD--PSI scheme \eqref{eq:ptd_lf_K}--\eqref{eq:ptd_lf_L}. If $\nu\le \frac13$, then
\[
\|U^{n+1}\|_F\le \|U^n\|_F,
\qquad \text{for all } n\ge 0.
\]
\end{theorem}

The proof follows the same overall strategy as in the DtP case, with the reduced estimates arranged according to the PtD substeps. We first state a lemma that will be used in the proof of the theorem.

\begin{lemma}
\label{lem:ptd_metric}
Let $\nu\le \frac13$, and let $T\in\mathbb R^{r\times r}$ be symmetric and satisfy
\[
0\preceq T\preceq 4I.
\]
Suppose that $W\in\mathbb R^{r\times r}$ is symmetric and satisfies
\[
W\succeq I,
\qquad
W-I\preceq \nu^2\Bigl(T-\frac14T^2\Bigr).
\]
Then
\[
W^2\preceq I+\nu(1-\nu)T.
\]
\end{lemma}

The proof is deferred to Appendix~\ref{app:ptd_metric_proof}.

\begin{proof}
We write
\[
|U|_{1,x}^2:=\langle U,M_2U\rangle_F.
\]

\paragraph{Step 1.}
The PtD K-step has the same full-matrix form as the DtP K-step and therefore produces the same mixed $L^2$--$H^1_x$ estimate.

Let
\[
G_0:=I-\frac{\nu}{2}M_2,
\qquad
G_1:=\frac{\Delta t}{2\Delta x}M_1,
\qquad
P_{V^n}:=V^n(V^n)^T.
\]
Recall that
\[
\widetilde A=(V^n)^TAV^n.
\]
Since
\[
(U^nV^n)\widetilde A(V^n)^T
=
U^nV^n(V^n)^TAV^n(V^n)^T
=
U^nP_{V^n}AP_{V^n}
=
U^nAP_{V^n},
\]
and
\[
\lambda_{\max}(U^nV^n)(V^n)^T
=
\lambda_{\max}U^nP_{V^n}
=
\lambda_{\max}U^n,
\]
equation \eqref{eq:ptd_lf_full_1} can be written as
\[
U^{(1)}=(G_0U^n-G_1U^nA)P_{V^n}.
\]
This is exactly the same algebraic form as in the DtP proof. Therefore the argument of Step~1 in the proof of Theorem~\ref{thm:dtp_general} applies without change and gives
\begin{equation}
\label{eq:ptd_mixed_1}
\|U^{(1)}\|_F^2+\nu(1-\nu)|U^{(1)}|_{1,x}^2
\le
\|U^n\|_F^2.
\end{equation}

\paragraph{Step 2.}
After the K-step, write
\[
U^{(1)}=X^{n+1}S^{(1)}(V^n)^T,
\qquad
\widehat A:=\frac{1}{\lambda_{\max}}\widetilde A.
\]
Since $A$ is symmetric and $\lambda_{\max}=\|A\|_2$, we have
\[
\widehat A^T=\widehat A,
\qquad
\widehat A^2\preceq I.
\]
Define
\[
M_{1,X}:=(X^{n+1})^TM_1X^{n+1},
\qquad
M_{2,X}:=(X^{n+1})^TM_2X^{n+1},
\]
and
\[
W_X:=I+\frac{\nu^2}{4}M_{1,X}^TM_{1,X}.
\]
Then
\[
M_{1,X}^T=-M_{1,X},
\qquad
0\preceq M_{2,X}\preceq 4I,
\qquad
W_X=I-\frac{\nu^2}{4}M_{1,X}^2\succeq I.
\]
Moreover,
\begin{equation}
\label{eq:ptd_h1_identity}
|U^{(1)}|_{1,x}^2
=
\operatorname{tr}\bigl((S^{(1)})^TM_{2,X}S^{(1)}\bigr).
\end{equation}
Indeed,
\[
|U^{(1)}|_{1,x}^2
=
\bigl\langle X^{n+1}S^{(1)}(V^n)^T,\,
M_2X^{n+1}S^{(1)}(V^n)^T\bigr\rangle_F
=
\operatorname{tr}\bigl((S^{(1)})^TM_{2,X}S^{(1)}\bigr).
\]

Define
\[
\mathcal C_X(L):=L-\frac{\nu}{2\lambda_{\max}}M_{1,X}LA,
\qquad L\in\mathbb R^{r\times N_v},
\]
and
\[
\mathcal C(Z):=Z+\frac{\nu}{2}M_{1,X}Z\widehat A,
\qquad Z\in\mathbb R^{r\times r}.
\]
Using \eqref{eq:ptd_lf_L}, we have
\[
L^n=S^{(2)}(V^n)^T,
\qquad
L^{n+1}=\mathcal C_X(L^n),
\qquad
U^{n+1}=X^{n+1}L^{n+1}.
\]
Since $X^{n+1}$ has orthonormal columns,
\[
\|U^{n+1}\|_F=\|L^{n+1}\|_F.
\]

For any $L\in\mathbb R^{r\times N_v}$, the matrix $L^TM_{1,X}L$ is skew-symmetric. Since $A$ is symmetric,
\[
\operatorname{tr}(L^TM_{1,X}LA)=0.
\]
Therefore,
\[
\|\mathcal C_X(L)\|_F^2
=
\|L\|_F^2+\frac{\nu^2}{4\lambda_{\max}^2}\|M_{1,X}LA\|_F^2.
\]
Using $A^2\preceq \lambda_{\max}^2I$, we obtain
\[
\|\mathcal C_X(L)\|_F^2
\le
\|L\|_F^2+\frac{\nu^2}{4}\|M_{1,X}L\|_F^2
=
\langle L,W_XL\rangle_F.
\]
Applying this to $L=L^n=S^{(2)}(V^n)^T$ and using $(V^n)^TV^n=I_r$, we find
\begin{equation}
\label{eq:ptd_step2_bound}
\|U^{n+1}\|_F^2
\le
\operatorname{tr}\bigl((S^{(2)})^TW_XS^{(2)}\bigr).
\end{equation}

Using \eqref{eq:ptd_lf_S}, we may write
\[
S^{(2)}=\mathcal C(S^{(1)}).
\]
Since $W_X=I-\frac{\nu^2}{4}M_{1,X}^2$ is a polynomial
in $M_{1,X}^2$, it commutes with $M_{1,X}$.
Consequently $W_XM_{1,X}$ is skew-symmetric,
because $W_X$ is symmetric, $M_{1,X}$ is skew-symmetric.
Hence the cross terms cancel when we expand
\[
\langle \mathcal C(Z),W_X\mathcal C(Z)\rangle_F.
\]
We obtain
\[
\langle \mathcal C(Z),W_X\mathcal C(Z)\rangle_F
=
\langle Z,W_XZ\rangle_F
+
\frac{\nu^2}{4}
\operatorname{tr}\bigl(\widehat A Z^TM_{1,X}^TW_XM_{1,X}Z\widehat A\bigr).
\]
The matrix
\[
Z^TM_{1,X}^TW_XM_{1,X}Z
\]
is positive semidefinite. Since $\widehat A^2\preceq I$, it follows that
\[
\operatorname{tr}\bigl(\widehat A Z^TM_{1,X}^TW_XM_{1,X}Z\widehat A\bigr)
\le
\operatorname{tr}\bigl(Z^TM_{1,X}^TW_XM_{1,X}Z\bigr).
\]
Moreover, since $M_{1,X}^T=-M_{1,X}$ and $W_X$ commutes
with $M_{1,X}$,
\begin{align*}
\frac{\nu^2}{4}M_{1,X}^TW_XM_{1,X}
&= -\frac{\nu^2}{4}M_{1,X}W_XM_{1,X}
 = -\frac{\nu^2}{4}M_{1,X}^2W_X
 = (W_X-I)W_X = W_X^2-W_X.
\end{align*}
Therefore,
\[
\langle \mathcal C(Z),W_X\mathcal C(Z)\rangle_F
\le
\langle Z,W_X^2Z\rangle_F.
\]
Applying this with $Z=S^{(1)}$ and using \eqref{eq:ptd_step2_bound}, we obtain
\begin{equation}
\label{eq:ptd_step2_final}
\|U^{n+1}\|_F^2
\le
\operatorname{tr}\bigl((S^{(1)})^TW_X^2S^{(1)}\bigr).
\end{equation}

\paragraph{Step 3.}
Let
\[
P_{X^{n+1}}:=X^{n+1}(X^{n+1})^T.
\]
Since $P_{X^{n+1}}\preceq I$, identity \eqref{eq:M1M2} gives
\[
M_{1,X}^TM_{1,X}
=
(X^{n+1})^TM_1^TP_{X^{n+1}}M_1X^{n+1}
\preceq
(X^{n+1})^T(4M_2-M_2^2)X^{n+1}.
\]
Moreover,
\[
M_{2,X}^2
=
(X^{n+1})^TM_2X^{n+1}(X^{n+1})^TM_2X^{n+1}
=
(X^{n+1})^TM_2P_{X^{n+1}}M_2X^{n+1}
\preceq
(X^{n+1})^TM_2^2X^{n+1},
\]
since $M_2\succeq 0$. Therefore,
\[
M_{1,X}^TM_{1,X}\preceq 4M_{2,X}-M_{2,X}^2.
\]
Set
\[
T:=M_{2,X}.
\]
Then
\[
W_X-I
=
\frac{\nu^2}{4}M_{1,X}^TM_{1,X}
\preceq
\nu^2\Bigl(T-\frac14T^2\Bigr).
\]
Since $W_X\succeq I$ and $0\preceq T\preceq 4I$, Lemma~\ref{lem:ptd_metric} yields
\[
W_X^2\preceq I+\nu(1-\nu)M_{2,X}.
\]
Combining this with \eqref{eq:ptd_step2_final} and \eqref{eq:ptd_h1_identity}, we obtain
\[
\|U^{n+1}\|_F^2
\le
\|U^{(1)}\|_F^2+\nu(1-\nu)|U^{(1)}|_{1,x}^2.
\]
Finally, \eqref{eq:ptd_mixed_1} gives
\[
\|U^{n+1}\|_F^2\le \|U^n\|_F^2.
\]
This proves the $L^2$-stability of the PtD--PSI scheme for $\nu\le \frac13$.
\end{proof}

\begin{remark}[A second-order PtD variant: single-mode calculation only] In the PtD formulation, we also provide a single-mode CFL result for the PSI scheme using the Strang splitting and SSP-RK2 in each substep.

Assume that $U^n=\bx_m\bv_k^T$ as in Proposition~\ref{prop:ptd_single_mode}, and keep the notation there. For the half-step quantity, let
\[
p_{m,k}^{(\frac12)}:=1-\frac{\nu}{2}y_m-i\frac{\nu_k}{2}z_m.
\]
Then one such step gives
\[
U^{n+1}=g^{\rm PtD,Strang}_{m,k}U^n,
\]
with amplification factor
\[
g^{\rm PtD,Strang}_{m,k}
=
R_{\mathrm{SSP2}}\!\bigl(p_{m,k}^{(\frac12)}\bigr)\,
R_{\mathrm{SSP2}}\!\Bigl(1+i\frac{\nu_k}{2}z_m\Bigr)\,
R_{\mathrm{SSP2}}(1-i\nu_k z_m)\,
R_{\mathrm{SSP2}}\!\Bigl(1+i\frac{\nu_k}{2}z_m\Bigr)\,
R_{\mathrm{SSP2}}\!\bigl(p_{m,k}^{(\frac12)}\bigr).
\]
Equivalently,
\[
|g^{\rm PtD,Strang}_{m,k}|^2
=
\bigl|R_{\mathrm{SSP2}}(p_{m,k}^{(\frac12)})\bigr|^4
\bigl|R_{\mathrm{SSP2}}(1+i\frac{\nu_k}{2}z_m)\bigr|^4
\bigl|R_{\mathrm{SSP2}}(1-i\nu_k z_m)\bigr|^2.
\]

If, in addition, $|\lambda_k|=\lambda_{\max}$, then $|\nu_k|=\nu$, and the above expression becomes a scalar function of $y_m\in[0,2]$ and $\nu$. A direct numerical scan of this function indicates that
\[
|g^{\rm PtD,Strang}_{m,k}|\le 1
\qquad\text{for all } y_m\in[0,2]
\]
whenever
\[
\nu \lesssim 2.
\]
\end{remark}

\section{Parabolic problem}
\label{sec:parabolic}
\label{sec:parab}

We now study the parabolic semi-discretization \eqref{eq:parabolic_semi_discrete_x_revised3}. It is well-known that this full-tensor scheme coupled with the backward Euler is unconditionally $L^2$-stable. We ask whether the corresponding low-rank PSI schemes can retain this property. The key issue is the time discretization used in the $S$-step. To identify the appropriate variant, we first consider a $\theta$-scheme for the $S$-step and use a single-mode calculation to show that the all-backward Euler variant is not unconditionally stable. This singles out the hybrid DtP scheme, with backward Euler in the $K$-step and $L$-step and forward Euler in the $S$-step. We then prove unconditional $L^2$-stability for this hybrid DtP scheme and finally transfer the same conclusion to the PtD formulation.

\subsection{Full-tensor backward Euler scheme}
\label{subsec:par_full_be}

Applying backward Euler to \eqref{eq:parabolic_semi_discrete_x_revised3} gives
\begin{equation}
\label{eq:par_full_be}
U^{n+1}=U^n-\tau M_2U^{n+1}A,
\qquad
\tau:=\frac{\Delta t}{(\Delta x)^2}.
\end{equation}
The next theorem records the classical contractivity property that serves as the full-tensor baseline for the low-rank PSI schemes considered below.

\begin{theorem}
\label{thm:par_full_be}
The full-tensor backward Euler scheme \eqref{eq:par_full_be} is unconditionally $L^2$-stable:
\[
\|U^{n+1}\|_F\le \|U^n\|_F
\qquad
\text{for all }\Delta t>0.
\]
\end{theorem}

\begin{proof}
Rearranging \eqref{eq:par_full_be}, we have
\[
U^n=U^{n+1}+\tau M_2U^{n+1}A.
\]
Taking the Frobenius norm squared gives
\[
\|U^n\|_F^2
=
\|U^{n+1}\|_F^2
+2\tau \langle U^{n+1},M_2U^{n+1}A\rangle_F
+\tau^2\|M_2U^{n+1}A\|_F^2.
\]
Since $M_2\succeq 0$ and $A\succeq 0$, we may write
\[
\langle U^{n+1},M_2U^{n+1}A\rangle_F
=
\|M_2^{1/2}U^{n+1}A^{1/2}\|_F^2
\ge 0.
\]
Therefore
\[
\|U^{n+1}\|_F\le \|U^n\|_F.
\]
This holds for all $\Delta t>0$.
\end{proof}

\subsection{DtP formulation}
\label{subsec:par_dtp}

\subsubsection{A $\theta$-scheme for the $S$-step and single-mode analysis}
\label{subsec:par_theta}

We begin with the DtP formulation for the parabolic semi-discretization \eqref{eq:parabolic_semi_discrete_x_revised3}. To compare the forward and backward treatments of the $S$-step, we introduce a standard $\theta$-scheme. Thus we consider
\[
\dot U=-\frac{1}{(\Delta x)^2}M_2UA,
\qquad
A=A^T\succeq 0.
\]
Let
\[
U(t)=X(t)S(t)V(t)^T,
\qquad
X\in\mathbb R^{N_x\times r},
\quad
S\in\mathbb R^{r\times r},
\quad
V\in\mathbb R^{N_v\times r},
\]
with
\[
X^TX=I_r,
\qquad
V^TV=I_r.
\]
Set
\[
K^n:=X^nS^n,
\qquad
\widetilde A:=(V^n)^TAV^n,
\qquad
\tau:=\frac{\Delta t}{(\Delta x)^2}.
\]
For a parameter $\theta\in[0,1]$, we keep backward Euler in the $K$-step and $L$-step, and let the $S$-step interpolate between forward and backward Euler:
\begin{equation}
\label{eq:par_theta_K}
K^{n+1}=K^n-\tau M_2K^{n+1}\widetilde A,
\qquad
K^{n+1}=X^{n+1}S^{(1)},
\end{equation}
\begin{equation}
\label{eq:par_theta_S}
S^{(2)}
=
S^{(1)}+\tau M_{2,X}\bigl((1-\theta)S^{(1)}+\theta S^{(2)}\bigr)\widetilde A,
\qquad
M_{2,X}:=(X^{n+1})^TM_2X^{n+1},
\end{equation}
\begin{equation}
\label{eq:par_theta_L}
L^{n+1}=L^n-\tau M_{2,X}L^{n+1}A,
\qquad
L^n:=S^{(2)}(V^n)^T,
\qquad
(L^{n+1})^T=V^{n+1}(S^{n+1})^T.
\end{equation}
The value $\theta=1$ corresponds to the all-backward Euler variant, whereas $\theta=0$ gives the hybrid DtP scheme, with backward Euler in the $K$-step and $L$-step and forward Euler in the $S$-step.

To see which choice can yield unconditional stability, we examine the action of \eqref{eq:par_theta_K}--\eqref{eq:par_theta_L} on a single Fourier mode. We use the Fourier modes $\bx_m$ introduced in Section~\ref{subsec:dtp_single_mode}; here we only need
\[
M_2\bx_m=2y_m\bx_m,
\qquad
y_m=1-\cos\alpha_m\in[0,2].
\]
For an eigenvector $\bv_k$ of $A$ with
\[
A\bv_k=\lambda_k\bv_k,
\qquad
\lambda_k\ge 0,
\qquad
\|\bv_k\|_2=1,
\]
define
\[
\mu:=\lambda_{\max}\frac{\Delta t}{(\Delta x)^2},
\qquad
\mu_k:=\lambda_k\frac{\Delta t}{(\Delta x)^2},
\qquad
\psi_{m,k}:=2\mu_k y_m.
\]

\begin{proposition}
\label{prop:par_theta_single_mode}
Assume that $U^n=\bx_m\bv_k^T$. Then one step of the DtP $\theta$-scheme defined by
\eqref{eq:par_theta_K}--\eqref{eq:par_theta_L} gives
\[
U^{n+1}=g_\theta(\psi_{m,k})\,U^n,
\]
where
\[
g_\theta(\psi)
=
\frac{1+(1-\theta)\psi}{(1+\psi)^2(1-\theta\psi)}.
\]
In particular, if $\theta=0$, then
\[
g_0(\psi)=\frac{1}{1+\psi},
\qquad \psi\ge 0.
\]
If $\theta>0$ and $\theta \neq 1/2$, then $g_\theta(\psi)$ becomes unbounded as $\psi$ approaches $1/\theta$ from below.
\end{proposition}

The proof is given in Appendix~\ref{app:par_theta_single_mode_calc}. For $\theta=1$, we have
\[
g_1(\psi)=\frac{1}{(1+\psi)^2(1-\psi)},
\]
which has a singularity at $\psi=1$. Since
\[
\psi_{m,k}=2\mu_k y_m=4\mu_k\sin^2\Bigl(\frac{\alpha_m}{2}\Bigr),
\]
this singularity can be approached on sufficiently fine grids for suitable values of $\mu_k$. Thus the all-backward Euler variant is not unconditionally stable.

By contrast, for $\theta=0$ we have $g_0(\psi)=(1+\psi)^{-1} \le 1$ for $\psi\ge 0$. This singles out the hybrid DtP scheme for the analysis below. 

\subsubsection{The hybrid DtP scheme and its unconditional $L^2$-stability}
\label{subsec:par_dtp_hybrid}

We now fix $\theta=0$ in the $\theta$-scheme \eqref{eq:par_theta_K}--\eqref{eq:par_theta_L}. This gives the hybrid DtP scheme, with backward Euler in the $K$-step and $L$-step and forward Euler in the $S$-step. We note that a similar time discretization is already used in
\cite{bachmayr2021existence,coughlin2022efficient} but without a stability proof.

Let
\[
K^n:=X^nS^n,
\qquad
\widetilde A:=(V^n)^TAV^n,
\qquad
\tau:=\frac{\Delta t}{(\Delta x)^2}.
\]
Then the DtP hybrid update is
\begin{equation}
\label{eq:par_hybrid_dtp_K}
K^{n+1}=K^n-\tau M_2K^{n+1}\widetilde A,
\qquad
K^{n+1}=X^{n+1}S^{(1)},
\end{equation}
\begin{equation}
\label{eq:par_hybrid_dtp_S}
S^{(2)}=S^{(1)}+\tau M_{2,X}S^{(1)}\widetilde A,
\qquad
M_{2,X}:=(X^{n+1})^TM_2X^{n+1},
\end{equation}
\begin{equation}
\label{eq:par_hybrid_dtp_L}
L^{n+1}=L^n-\tau M_{2,X}L^{n+1}A,
\qquad
L^n:=S^{(2)}(V^n)^T,
\qquad
(L^{n+1})^T=V^{n+1}(S^{n+1})^T.
\end{equation}

The next theorem shows that the hybrid DtP scheme is unconditionally $L^2$-stable.
\begin{theorem}
\label{thm:par_hybrid_dtp}
The hybrid DtP scheme \eqref{eq:par_hybrid_dtp_K}--\eqref{eq:par_hybrid_dtp_L} is unconditionally $L^2$-stable:
\[
\|U^{n+1}\|_F\le \|U^n\|_F
\qquad
\text{for all }\Delta t>0.
\]
\end{theorem}

\begin{proof}
Substituting $K^{n+1}=X^{n+1}S^{(1)}$ and $K^n=X^nS^n$ into \eqref{eq:par_hybrid_dtp_K}, we obtain
\[
X^{n+1}S^{(1)}=X^nS^n-\tau M_2X^{n+1}S^{(1)}\widetilde A.
\]
Left-multiplying by $(X^{n+1})^T$ gives
\[
S^{(1)}=RS^n-\tau M_{2,X}S^{(1)}\widetilde A,
\qquad
R:=(X^{n+1})^TX^n.
\]
Combining this with \eqref{eq:par_hybrid_dtp_S}, we find
\[
S^{(2)}
=
\bigl(RS^n-\tau M_{2,X}S^{(1)}\widetilde A\bigr)
+\tau M_{2,X}S^{(1)}\widetilde A
=
RS^n.
\]
Let
\[
P_{X^{n+1}}:=X^{n+1}(X^{n+1})^T.
\]
Hence
\[
X^{n+1}S^{(2)}
=
X^{n+1}(X^{n+1})^TX^nS^n
=
P_{X^{n+1}}K^n.
\]
Since $X^{n+1}$ has orthonormal columns,
\[
\|S^{(2)}\|_F
=
\|X^{n+1}S^{(2)}\|_F
=
\|P_{X^{n+1}}K^n\|_F
\le
\|K^n\|_F
=
\|S^n\|_F
=
\|U^n\|_F.
\]

We next consider the $L$-step. Rearranging \eqref{eq:par_hybrid_dtp_L}, we have
\[
L^n=L^{n+1}+\tau M_{2,X}L^{n+1}A.
\]
Taking the Frobenius norm squared gives
\[
\|L^n\|_F^2
=
\|L^{n+1}\|_F^2
+2\tau \langle L^{n+1},M_{2,X}L^{n+1}A\rangle_F
+\tau^2\|M_{2,X}L^{n+1}A\|_F^2.
\]
Since $M_{2,X}\succeq 0$ and $A\succeq 0$,
\[
\langle L^{n+1},M_{2,X}L^{n+1}A\rangle_F
=
\operatorname{tr}\bigl((L^{n+1})^TM_{2,X}L^{n+1}A\bigr)
=
\|M_{2,X}^{1/2}L^{n+1}A^{1/2}\|_F^2
\ge 0.
\]
Therefore
\[
\|L^{n+1}\|_F\le \|L^n\|_F.
\]

Finally, since
\[
L^n=S^{(2)}(V^n)^T,
\qquad
U^{n+1}=X^{n+1}L^{n+1},
\]
and both $V^n$ and $X^{n+1}$ have orthonormal columns, we have
\[
\|L^n\|_F=\|S^{(2)}\|_F,
\qquad
\|U^{n+1}\|_F=\|L^{n+1}\|_F.
\]
Combining the above estimates yields
\[
\|U^{n+1}\|_F
=
\|L^{n+1}\|_F
\le
\|L^n\|_F
=
\|S^{(2)}\|_F
\le
\|U^n\|_F.
\]
This proves the stated $L^2$-stability result.
\end{proof}

\subsection{PtD formulation}
\label{subsec:par_ptd}

We finally turn to the PtD formulation for the parabolic problem. Starting from the semi-discrete system \eqref{eq:parabolic_semi_discrete_v}, we first apply the low-rank projection and then discretize the resulting $K$-, $S$-, and $L$-steps in $x$ by the same central difference operator as before. For the present linear constant-coefficient parabolic model, the low-rank projection and the finite-difference discretization commute. Therefore, the stability result established in the previous subsection for the DtP formulation applies here {\it verbatim}. This is in contrast to the hyperbolic case, where different spatial discretizations are used in different substeps.

To see this, we start from the transposed form of \eqref{eq:parabolic_semi_discrete_v},
\[
\partial_t \bu^T(t,x)=\partial_{xx}\bu^T(t,x)A,
\]
and consider the rank-$r$ approximation
\[
\bu^T(t,x)=\bx^T(t,x)S(t)V(t)^T,
\]
where $\bx(t,x)\in\mathbb R^r$, $S(t)\in\mathbb R^{r\times r}$, and $V(t)\in\mathbb R^{N_v\times r}$ satisfy
\[
\langle \bx\bx^T\rangle_x=I_r,
\qquad
V^TV=I_r.
\]
Let
\[
\bk^T(t,x):=\bx^T(t,x)S(t),
\qquad
L(t):=S(t)V(t)^T,
\qquad
\widetilde A:=V^TAV.
\]
Then the projector-splitting subproblems read
\begin{align}
\partial_t \bk^T &= \partial_{xx}\bk^T\,\widetilde A,
\label{eq:par_ptd_cont_K}
\\
\dot S &= -\langle \bx\,\partial_{xx}\bx^T\rangle_x\,S\widetilde A,
\label{eq:par_ptd_cont_S}
\\
\dot L &= \langle \bx\,\partial_{xx}\bx^T\rangle_x\,LA.
\label{eq:par_ptd_cont_L}
\end{align}

Discretizing these equations in $x$ by the second-order central difference scheme, and introducing
\[
X(t):=\sqrt{\Delta x}
\begin{bmatrix}
\bx_1(t)^T\\
\vdots\\
\bx_{N_x}(t)^T
\end{bmatrix}
\in\mathbb R^{N_x\times r},
\qquad
K(t):=X(t)S(t),
\qquad
M_{2,X}:=X^TM_2X,
\]
we obtain
\begin{align}
\dot K &= -\frac{1}{(\Delta x)^2}M_2K\widetilde A,
\label{eq:par_ptd_sd_K}
\\
\dot S &= \frac{1}{(\Delta x)^2}M_{2,X}S\widetilde A,
\label{eq:par_ptd_sd_S}
\\
\dot L &= -\frac{1}{(\Delta x)^2}M_{2,X}LA.
\label{eq:par_ptd_sd_L}
\end{align}
These are exactly the same algebraic subproblems as in the DtP formulation.

Applying backward Euler in the $K$-step and $L$-step and forward Euler in the $S$-step therefore gives the hybrid PtD scheme
\begin{equation}
\label{eq:par_hybrid_ptd_K}
K^{n+1}=K^n-\tau M_2K^{n+1}\widetilde A,
\qquad
K^{n+1}=X^{n+1}S^{(1)},
\end{equation}
\begin{equation}
\label{eq:par_hybrid_ptd_S}
S^{(2)}=S^{(1)}+\tau M_{2,X}S^{(1)}\widetilde A,
\qquad
M_{2,X}:=(X^{n+1})^TM_2X^{n+1},
\end{equation}
\begin{equation}
\label{eq:par_hybrid_ptd_L}
L^{n+1}=L^n-\tau M_{2,X}L^{n+1}A,
\qquad
L^n:=S^{(2)}(V^n)^T,
\qquad
(L^{n+1})^T=V^{n+1}(S^{n+1})^T,
\end{equation}
where
\[
\tau:=\frac{\Delta t}{(\Delta x)^2},
\qquad
\widetilde A:=(V^n)^TAV^n.
\]

\begin{theorem}
\label{thm:par_hybrid_ptd}
The hybrid PtD scheme \eqref{eq:par_hybrid_ptd_K}--\eqref{eq:par_hybrid_ptd_L} is unconditionally $L^2$-stable:
\[
\|U^{n+1}\|_F\le \|U^n\|_F
\qquad
\text{for all }\Delta t>0.
\]
\end{theorem}

\begin{proof}
Since \eqref{eq:par_ptd_sd_K}--\eqref{eq:par_ptd_sd_L} are algebraically identical to the corresponding DtP subproblems, the hybrid PtD scheme \eqref{eq:par_hybrid_ptd_K}--\eqref{eq:par_hybrid_ptd_L} has exactly the same form as the hybrid DtP scheme \eqref{eq:par_hybrid_dtp_K}--\eqref{eq:par_hybrid_dtp_L}. The conclusion therefore follows directly from Theorem~\ref{thm:par_hybrid_dtp}.
\end{proof}

\section{Conclusion}
\label{sec:conclusion}

We have analyzed the stability of low-rank PSI schemes for linear hyperbolic and parabolic prototype equations.

For the hyperbolic problem, the full-tensor LF scheme with forward Euler satisfies the classical hyperbolic CFL condition $\nu\le 1$. For the corresponding low-rank PSI schemes, single-mode calculations lead to the bound $\nu\le 1/3$, and we proved $L^2$-stability under this bound in both the DtP and PtD formulations.

For the parabolic problem, while the full-tensor centered finite difference scheme with backward Euler is unconditionally $L^2$-stable, we proved that the corresponding PSI scheme is not. However, a hybrid PSI using backward Euler in the $K$- and $L$-steps and forward Euler in the $S$-step is unconditionally $L^2$-stable in both the DtP and PtD formulations.

Although the analysis is carried out for prototype equations, the conclusions obtained here are informative for more general kinetic equations. Future directions include studying the stability of high-order splitting schemes, high-order temporal and spatial discretizations.

\section*{Acknowledgement}
This work was partially supported by DOE grant DE-SC0023164, NSF grants DMS-2409858 and IIS-2433957, and DoD MURI grant FA9550-24-1-0254.

\section*{Data availability statement}
There is no data generated from this work.

\bibliographystyle{siamplain}
\bibliography{hu_bibtex}

@article{einkemmer2025review,
  title={A review of low-rank methods for time-dependent kinetic simulations},
  author={Einkemmer, Lukas and Kormann, Katharina and Kusch, Jonas and McClarren, Ryan G and Qiu, Jing-Mei},
  journal={Journal of Computational Physics},
  volume={538},
  pages={114191},
  year={2025},
  publisher={Elsevier}
}

@article{EHZ25,
	author = {L. Einkemmer and J. Hu and S. Zhang},
	journal = {J. Comput. Phys.},
	pages = {114112},
	title = {Asymptotic-preserving dynamical low-rank method for the stiff nonlinear {B}oltzmann equation},
	volume = {538},
	year = {2025}}

@article{CKL24,
	author = {G. Ceruti and J. Kusch and C. Lubich},
	journal = {SIAM J. Sci. Comput.},
	pages = {B205--B228},
	title = {A parallel rank-adaptive integrator for dynamical low-rank approximation},
	volume = {46},
	year = {2024}}

@article{CKL22,
	author = {G. Ceruti and J. Kusch and C. Lubich},
	journal = {BIT Numer. Math.},
	pages = {1149--1174},
	title = {A rank-adaptive robust integrator for dynamical low-rank approximation},
	volume = {62},
	year = {2022}}

@article{KEC23,
	author = {J. Kusch and L. Einkemmer and G. Ceruti},
	journal = {SIAM J. Sci. Comput.},
	pages = {A1--A24},
	title = {On the stability of robust dynamical low-rank approximations for hyperbolic problems},
	volume = {45},
	year = {2023}}

@article{EL18,
	author = {L. Einkemmer and C. Lubich},
	journal = {SIAM J. Sci. Comput.},
	pages = {B1330--B1360},
	title = {A low-rank projector-splitting integrator for the {V}lasov-{P}oisson equation},
	volume = {40},
	year = {2018}}

@article{EHW21,
	author = {L. Einkemmer and J. Hu and Y. Wang},
	journal = {J. Comput. Phys.},
	pages = {110353},
	title = {{An asymptotic-preserving dynamical low-rank method for the multi-scale multi-dimensional linear transport equation}},
	volume = {439},
	year = {2021}}

@article{LO14,
	author = {C. Lubich and I. Oseledets},
	journal = {BIT Numer. Math.},
	pages = {171--188},
	title = {A projector-splitting integrator for dynamical low-rank approximation},
	volume = {54},
	year = {2014}}

@book{GKS11,
	author = {S. Gottlieb and D. Ketcheson and C.-W. Shu},
	publisher = {World Scientific},
	title = {Strong Stability Preserving Runge-Kutta and Multistep Time Discretizations},
	year = {2011}}

@incollection{Villani02,
	author = {C. Villani},
	booktitle = {Handbook of Mathematical Fluid Mechanics},
	editor = {S. Friedlander and D. Serre},
	pages = {71--305},
	publisher = {North-Holland},
	title = {A review of mathematical topics in collisional kinetic theory},
	volume = {I},
	year = {2002}}

@article{kazashi2021stability,
  title={Stability properties of a projector-splitting scheme for dynamical low rank approximation of random parabolic equations},
  author={Kazashi, Yoshihito and Nobile, Fabio and Vidli{\v{c}}kov{\'a}, Eva},
  journal={Numerische Mathematik},
  volume={149},
  number={4},
  pages={973--1024},
  year={2021},
  publisher={Springer}
}

@article{bachmayr2021existence,
  title={Existence of dynamical low-rank approximations to parabolic problems},
  author={Bachmayr, Markus and Eisenmann, Henrik and Kieri, Emil and Uschmajew, Andr{\'e}},
  journal={Mathematics of Computation},
  volume={90},
  number={330},
  pages={1799--1830},
  year={2021}
}

@article{coughlin2022efficient,
  title={{Efficient dynamical low-rank approximation for the Vlasov-Amp{\`e}re-Fokker-Planck system}},
  author={Coughlin, Jack and Hu, Jingwei},
  journal={Journal of Computational Physics},
  volume={470},
  pages={111590},
  year={2022},
  publisher={Elsevier}
}

\begin{appendices}

\section{DtP formulation for the hyperbolic problem}
\label{app:dtp_hyper}

This appendix collects the details used in Section~\ref{subsec:hyper_dtp}.

\subsection{Proof of Proposition~\ref{prop:dtp_single_mode}}
\label{app:dtp_single_mode_calc}

\begin{proof}
Let $U^n=\bx_m\bv_k^T$, where $A\bv_k=\lambda_k\bv_k$ and $\|\bv_k\|_2=1$. Then
\[
P_{V^n}=\bv_k\bv_k^T,
\qquad
M_1\bx_m=2iz_m\bx_m,
\qquad
M_2\bx_m=2y_m\bx_m.
\]
Using \eqref{eq:dtp_lf_full_1}, we obtain
\[
U^{(1)}
=
\bx_m\bv_k^T
-\frac{\Delta t}{2\Delta x}\bigl(2iz_m\bx_m\,\lambda_k\bv_k^T+2\lambda_{\max}y_m\bx_m\bv_k^T\bigr)\bv_k\bv_k^T
=
p_{m,k}U^n,
\]
with
\[
p_{m,k}=1-\nu y_m-i\nu_k z_m,
\qquad
\nu_k=\lambda_k\frac{\Delta t}{\Delta x}.
\]
Since $U^{(1)}$ is a scalar multiple of $U^n$, it lies in the same rank-one subspace, and therefore
\[
P_{X^{n+1}}U^{(1)}=U^{(1)},
\qquad
U^{(1)}P_{V^n}=U^{(1)}.
\]
Using \eqref{eq:dtp_lf_full_2}, we find
\[
U^{(2)}
=
U^{(1)}+\frac{\Delta t}{2\Delta x}\bigl(2iz_m\,\lambda_k+2\lambda_{\max}y_m\bigr)U^{(1)}
=
(2-p_{m,k})U^{(1)}.
\]
Applying \eqref{eq:dtp_lf_full_3} in the same way gives
\[
U^{n+1}
=
p_{m,k}U^{(2)}
=
p_{m,k}^2(2-p_{m,k})U^n.
\]
This proves the amplification formula.

Taking absolute values gives
\[
|g^{\rm DtP}_{m,k}|^2
=
|p_{m,k}|^4|2-p_{m,k}|^2
=
\bigl((1-\nu y_m)^2+\nu_k^2z_m^2\bigr)^2
\bigl((1+\nu y_m)^2+\nu_k^2z_m^2\bigr).
\]
If $|\lambda_k|=\lambda_{\max}$, then $|\nu_k|=\nu$, and the identity $z_m^2=y_m(2-y_m)$ yields
\[
|p_{m,k}|^2
=
(1-\nu y_m)^2+\nu^2z_m^2
=
1+2y_m\nu(\nu-1),
\]
and
\[
|2-p_{m,k}|^2
=
(1+\nu y_m)^2+\nu^2z_m^2
=
1+2y_m\nu(\nu+1).
\]
Hence
\[
|g^{\rm DtP}_{m,k}|^2
=
h_{\rm DtP}(y_m,\nu)
=
\bigl[1+2y_m\nu(\nu-1)\bigr]^2
\bigl[1+2y_m\nu(\nu+1)\bigr].
\]

To identify the modewise threshold, differentiate with respect to $y$:
\[
\frac{\partial h_{\rm DtP}}{\partial y}
=
2\nu\bigl[1+2y\nu(\nu-1)\bigr]\bigl[3\nu-1+6y\nu(\nu^2-1)\bigr].
\]
For $0\le \nu\le \frac13$, the first factor is nonnegative on $[0,2]$ because
\[
1+2y\nu(\nu-1)=(1-\nu y)^2+\nu^2y(2-y)\ge 0,
\]
and the second factor is nonpositive on $[0,2]$. Hence $h_{\rm DtP}(y,\nu)$ is nonincreasing on $[0,2]$, and since $h_{\rm DtP}(0,\nu)=1$, we obtain
\[
h_{\rm DtP}(y,\nu)\le 1,
\qquad y\in[0,2].
\]
If $\nu>\frac13$, then
\[
\frac{\partial h_{\rm DtP}}{\partial y}(0,\nu)=2\nu(3\nu-1)>0,
\]
so $h_{\rm DtP}(y,\nu)>1$ for all sufficiently small $y>0$.
\end{proof}

\subsection{Proof of Lemma~\ref{lem:dtp_metric}}
\label{app:dtp_metric_proof}

\begin{proof}
Choose any orthonormal basis of $\mathbb R^{r\times r}$ with respect to the Frobenius inner product, and let $E\in\mathbb R^{r^2\times r^2}$ be the matrix representation of the operator $\mathcal E$ in this basis. Then
\[
\|E\|_2=\|\mathcal E\|_{F\to F}\le 1,
\]
and the adjoint $\mathcal E^*$ is represented by $E^T$. Likewise, the operators
\[
\mathcal H=2\mathcal I-(\mathcal E+\mathcal E^*),
\qquad
\mathcal W=\mathcal I-\nu(1-\nu)\mathcal H,
\qquad
\mathcal Q=\mathcal I+\nu(\mathcal I-\mathcal E)
\]
are represented by the matrices
\[
H:=2I-(E+E^T),
\qquad
W:=I-\nu(1-\nu)H,
\qquad
Q:=I+\nu(I-E).
\]
Because the quadratic-form order is preserved under this orthonormal identification, it suffices to prove
\[
Q^TWQ\preceq I+\nu(1-\nu)H.
\]

Set
\[
Z:=I-E,
\qquad
D^{\mathrm L}:=I-EE^T,
\qquad
D^{\mathrm R}:=I-E^TE.
\]
Since $\|E\|_2\le 1$, we have
\[
-2I\preceq E+E^T\preceq 2I,
\]
and therefore
\[
0\preceq H=2I-(E+E^T)\preceq 4I.
\]
Hence, for $0<\nu\le \frac13$,
\[
W=I-\nu(1-\nu)H \succeq (1-4\nu(1-\nu))I \succeq \frac19 I,
\]
so $W$ is symmetric positive definite.

A direct expansion yields
\[
Z+Z^T=H,
\qquad
Z^TZ=H-D^{\mathrm R},
\qquad
ZZ^T=H-D^{\mathrm L},
\qquad
Z^TH+HZ=H^2+D^{\mathrm L}-D^{\mathrm R}.
\]
Using these identities, we compute
\begin{align*}
Q^TWQ
&=(I+\nu Z^T)(I-\nu(1-\nu)H)(I+\nu Z)
\\
&=I+2\nu^2H-\nu^2(1-\nu)H^2
-\nu^2(1-\nu)D^{\mathrm L}
-\nu^3D^{\mathrm R}
-\nu^3(1-\nu)Z^THZ.
\end{align*}
Because $D^{\mathrm L}\succeq 0$, $D^{\mathrm R}\succeq 0$, $H\succeq 0$, and thus $Z^THZ\succeq 0$, we obtain
\[
Q^TWQ
\preceq
I+2\nu^2H-\nu^2(1-\nu)H^2.
\]

It remains to show that
\[
I+2\nu^2H-\nu^2(1-\nu)H^2
\preceq
I+\nu(1-\nu)H.
\]
Both sides are polynomials in the symmetric matrix $H$.
Since $0\preceq H\preceq 4I$,
the spectrum of $H$ lies in $[0,4]$, and by the
spectral mapping theorem it suffices to verify the
scalar inequality
\[
1+2\nu^2y-\nu^2(1-\nu)y^2
\le
1+\nu(1-\nu)y,
\qquad y\in[0,4].
\]
This is equivalent to
\[
\nu y\bigl(3\nu-1-\nu(1-\nu)y\bigr)\le 0.
\]
Since $0<\nu\le \frac13$ and $y\ge 0$, the above inequality holds for all $y\in[0,4]$. Therefore
\[
Q^TWQ\preceq I+\nu(1-\nu)H.
\]
Transferring this back to the Frobenius-orthonormal basis proves the operator inequality
\[
\mathcal Q^*\mathcal W\mathcal Q
\preceq
\mathcal I+\nu(1-\nu)\mathcal H.
\]
\end{proof}

\section{PtD formulation for the hyperbolic problem}
\label{app:ptd_hyper}

This appendix collects the details used in Section~\ref{subsec:hyper_ptd}.

\subsection{Proof of Proposition~\ref{prop:ptd_single_mode}}
\label{app:ptd_single_mode_calc}

\begin{proof}
Let $U^n=\bx_m\bv_k^T$, where $A\bv_k=\lambda_k\bv_k$ and $\|\bv_k\|_2=1$. Then
\[
U^nV^n=\bx_m,
\qquad
\widetilde A=(V^n)^TAV^n=\lambda_k,
\qquad
M_1\bx_m=2iz_m\bx_m,
\qquad
M_2\bx_m=2y_m\bx_m.
\]
Using \eqref{eq:ptd_lf_full_1}, we obtain
\[
U^{(1)}
=
\bx_m\bv_k^T-\frac{\Delta t}{2\Delta x}\bigl(2iz_m\lambda_k+2\lambda_{\max}y_m\bigr)\bx_m\bv_k^T
=
p_{m,k}U^n,
\]
with
\[
p_{m,k}=1-\nu y_m-i\nu_k z_m,
\qquad
\nu_k=\lambda_k\frac{\Delta t}{\Delta x}.
\]
Since $U^{(1)}$ is again a scalar multiple of $\bx_m\bv_k^T$, the updated column space is still spanned by $\bx_m$. Hence $X^{n+1}(X^{n+1})^TM_1(U^{(1)}V^n)=M_1(U^{(1)}V^n)$, and \eqref{eq:ptd_lf_full_2} gives
\[
U^{(2)}
=
U^{(1)}+\frac{\Delta t}{2\Delta x}2iz_m\lambda_k\,U^{(1)}
=
(1+i\nu_k z_m)U^{(1)}.
\]
Using \eqref{eq:ptd_lf_full_3} in the same way, we obtain
\[
U^{n+1}
=
(1-i\nu_k z_m)U^{(2)}
=
p_{m,k}(1+i\nu_k z_m)(1-i\nu_k z_m)U^n.
\]
This proves the amplification formula.

Taking absolute values gives
\[
|g^{\rm PtD}_{m,k}|^2
=
|p_{m,k}|^2(1+\nu_k^2z_m^2)^2
=
\bigl((1-\nu y_m)^2+\nu_k^2z_m^2\bigr)\bigl(1+\nu_k^2z_m^2\bigr)^2.
\]
If $|\lambda_k|=\lambda_{\max}$, then $|\nu_k|=\nu$, and $z_m^2=y_m(2-y_m)$ gives
\[
|g^{\rm PtD}_{m,k}|^2
=
h_{\rm PtD}(y_m,\nu)
=
\bigl[1+2y_m\nu(\nu-1)\bigr]\bigl[1+\nu^2y_m(2-y_m)\bigr]^2.
\]

To identify the modewise threshold, differentiate with respect to $y$:
\[
\frac{\partial h_{\rm PtD}}{\partial y}
=
2\nu\bigl[1+\nu^2y(2-y)\bigr]g(y;\nu),
\]
where
\[
g(y;\nu):=(3\nu-1)+2\nu(3\nu^2-3\nu-1)y+5\nu^2(1-\nu)y^2.
\]
For $0\le \nu\le \frac13$, the polynomial $g(\cdot;\nu)$ is convex, and
\[
g(0;\nu)=3\nu-1\le 0,
\qquad
g(2;\nu)=-1-\nu+8\nu^2(1-\nu)\le -\frac19-\nu<0.
\]
Hence $g(y;\nu)\le 0$ for all $y\in[0,2]$, so $h_{\rm PtD}(y,\nu)$ is nonincreasing on $[0,2]$. Since $h_{\rm PtD}(0,\nu)=1$, we obtain
\[
h_{\rm PtD}(y,\nu)\le 1,
\qquad y\in[0,2].
\]
If $\nu>\frac13$, then
\[
\frac{\partial h_{\rm PtD}}{\partial y}(0,\nu)=2\nu(3\nu-1)>0,
\]
so $h_{\rm PtD}(y,\nu)>1$ for all sufficiently small $y>0$.
\end{proof}

\subsection{Proof of Lemma~\ref{lem:ptd_metric}}
\label{app:ptd_metric_proof}

\begin{proof}
Set
\[
K:=W-I.
\]
Then $K\succeq 0$ and, by assumption,
\[
0\preceq K\preceq \nu^2\Bigl(T-\frac14T^2\Bigr).
\]
Since $0\preceq T\preceq 4I$, the scalar bound
\[
0\le c-\frac14c^2\le 1,
\qquad c\in[0,4],
\]
implies
\[
0\preceq T-\frac14T^2\preceq I.
\]
Hence
\[
0\preceq K\preceq \nu^2I.
\]

Define, for $t\in[0,\nu^2]$,
\[
h_\nu(t):=2-2\sqrt{1-t/\nu^2}.
\]
Since $t\mapsto 1-t/\nu^2$ is affine and order-reversing,
and $s\mapsto\sqrt{s}$ is operator-monotone increasing on
$[0,+\infty)$, the composition $t\mapsto\sqrt{1-t/\nu^2}$
is operator-monotone decreasing on $[0,\nu^2]$, and
therefore $h_\nu(t)=2-2\sqrt{1-t/\nu^2}$ is
operator-monotone increasing on $[0,\nu^2]$ with respect
to the Loewner order. Applying functional calculus to
$0\preceq K\preceq \nu^2I$, we obtain
\[
h_\nu(K)
\preceq
h_\nu\!\left(\nu^2\Bigl(T-\frac14T^2\Bigr)\right).
\]
For a scalar $c\in[0,4]$,
\[
h_\nu\!\left(\nu^2\Bigl(c-\frac14c^2\Bigr)\right)
=
2-2\sqrt{1-c+\frac14c^2}
=
2-|2-c|
\le c.
\]
Therefore, by functional calculus for the symmetric matrix $T$,
\begin{equation}
\label{eq:ptd_h_bound}
h_\nu(K)\preceq T.
\end{equation}

We next prove the scalar inequality
\[
(1+t)^2\le 1+\nu(1-\nu)h_\nu(t),
\qquad 0\le t\le \nu^2.
\]
Write
\[
t=\nu^2(1-s^2),
\qquad s\in[0,1].
\]
Then $h_\nu(t)=2(1-s)$, and
\[
(1+t)^2-\bigl(1+\nu(1-\nu)h_\nu(t)\bigr)
=
\nu(1-s)\Bigl(2\nu(1+s)+\nu^3(1-s)(1+s)^2-2(1-\nu)\Bigr).
\]
Set
\[
f_\nu(s):=2\nu(1+s)+\nu^3(1-s)(1+s)^2.
\]
Then
\[
f_\nu'(s)=2\nu+\nu^3(1+s)(1-3s)\ge 2\nu-4\nu^3>0
\qquad \text{for } 0\le s\le 1,
\]
whenever $0<\nu\le \frac13$. Thus $f_\nu$ is increasing on $[0,1]$, and hence
\[
f_\nu(s)\le f_\nu(1)=4\nu\le 2(1-\nu)
\qquad \text{for } 0<\nu\le \frac13.
\]
Therefore
\[
(1+t)^2\le 1+\nu(1-\nu)h_\nu(t),
\qquad 0\le t\le \nu^2.
\]

Applying functional calculus to the symmetric matrix $K$, we obtain
\[
W^2=(I+K)^2\preceq I+\nu(1-\nu)h_\nu(K).
\]
Using \eqref{eq:ptd_h_bound}, we conclude that
\[
W^2\preceq I+\nu(1-\nu)T.
\]
\end{proof}

\section{A $\theta$-scheme for the DtP parabolic problem}
\label{app:par_theta}

\subsection{Proof of Proposition~\ref{prop:par_theta_single_mode}}
\label{app:par_theta_single_mode_calc}

\begin{proof}
Let $U^n=\bx_m\bv_k^T$, where $A\bv_k=\lambda_k\bv_k$, $\lambda_k\ge 0$, and $\|\bv_k\|_2=1$. Then
\[
\widetilde A=(V^n)^TAV^n=\lambda_k,
\qquad
M_2\bx_m=2y_m\bx_m.
\]
Set
\[
\psi_{m,k}:=2\mu_k y_m,
\qquad
\mu_k:=\lambda_k\frac{\Delta t}{(\Delta x)^2}.
\]

For the K-step \eqref{eq:par_theta_K}, write
\[
K^{n+1}=\gamma_1\bx_m.
\]
Then
\[
\gamma_1\bx_m
=
\bx_m-\tau M_2(\gamma_1\bx_m)\lambda_k
=
\bx_m-\psi_{m,k}\gamma_1\bx_m,
\]
so
\[
\gamma_1=\frac{1}{1+\psi_{m,k}}.
\]
Since $\psi_{m,k}\ge 0$, we may take
\[
X^{n+1}=\bx_m,
\qquad
S^{(1)}=\frac{1}{1+\psi_{m,k}}.
\]
Moreover,
\[
M_{2,X}=(X^{n+1})^TM_2X^{n+1}=2y_m.
\]

For the S-step \eqref{eq:par_theta_S}, we obtain
\[
S^{(2)}
=
S^{(1)}+\psi_{m,k}\bigl((1-\theta)S^{(1)}+\theta S^{(2)}\bigr).
\]
Hence
\[
(1-\theta\psi_{m,k})S^{(2)}
=
\bigl(1+(1-\theta)\psi_{m,k}\bigr)S^{(1)},
\]
and therefore
\[
S^{(2)}
=
\frac{1+(1-\theta)\psi_{m,k}}{1-\theta\psi_{m,k}}\,S^{(1)}.
\]

Finally, for the L-step \eqref{eq:par_theta_L}, let
\[
L^{n+1}=\gamma_3\bv_k^T.
\]
Since $L^n=S^{(2)}\bv_k^T$ and $\bv_k^TA=\lambda_k\bv_k^T$, we have
\[
\gamma_3\bv_k^T
=
S^{(2)}\bv_k^T-\tau (2y_m)\gamma_3\bv_k^TA
=
S^{(2)}\bv_k^T-\psi_{m,k}\gamma_3\bv_k^T.
\]
Thus
\[
\gamma_3=\frac{S^{(2)}}{1+\psi_{m,k}}.
\]
Combining the three substeps gives
\[
U^{n+1}
=
\frac{1+(1-\theta)\psi_{m,k}}{(1+\psi_{m,k})^2(1-\theta\psi_{m,k})}\,U^n.
\]
This is precisely the stated amplification factor.
\end{proof}

\end{appendices}

\end{document}